\documentclass[11pt]{amsart}

\usepackage{geometry}                % See geometry.pdf to learn the layout options. There are lots.
\usepackage{graphicx}
\usepackage{framed}

\usepackage{algorithm}
\usepackage{algorithmicx}
\usepackage{algpseudocode}
\usepackage{multicol}

\usepackage{epstopdf}
\usepackage{amsmath,amssymb,amsthm,amsfonts}

\usepackage{mathrsfs}
\usepackage{hyperref}
\usepackage{tikz}
\usetikzlibrary{arrows,positioning,shapes}
\usepackage{enumitem}

\usepackage{verbatim}

\usepackage{array}

\theoremstyle{plain}  % Text in italic, and adds extra space before and after
\newtheorem{thm}{Theorem}[section]

\newtheorem{prop}[thm]{Proposition}

\theoremstyle{definition}  %  Text in roman, and adds extra space before and after
\newtheorem{defn}[thm]{Definition}

\theoremstyle{remark}  % Text in roman, but does not add extra space before or after
\newtheorem{rem}[thm]{Remark}

\newcommand{\setof}[1]{\left\{ {#1}\right\}}

\newcommand{\R}{{\mathbb{R}}}

\newcommand{\cP}{{\mathcal P}}

\newcommand{\cS}{{\mathcal S}}

\newcommand{\cX}{{\mathcal X}}

\newcommand{\cl}{\mathop{\mathrm{cl}}\nolimits}

\definecolor{gray85}{gray}{0.85} % 15%
\definecolor{gray8}{gray}{0.8} % 20%
\definecolor{gray7}{gray}{0.7} % 30%
\definecolor{gray6}{gray}{0.6} % 40%
\definecolor{gray5}{gray}{0.5} % 50%
\definecolor{gray4}{gray}{0.4} % 60%
\definecolor{gray35}{gray}{0.35} % 65%

\def\corcommstyle{\bf\small\tt}

\def\corrl #1<<#2||#3>>{
\if\visiblecomments y
  \begin{quote} {\corcommstyle $<<$COMMENT$>>$ {\color{red}#1\marginpar{!!}}\\$<<$OLD$<<$} \end{quote}

{\color{red} 
 #2
 }

  \begin{quote} {\corcommstyle ==NEW== } \end{quote}
   \noindent\hrulefill
 
\vspace{-10pt} 
 
 \noindent\hrulefill
 
 \vspace{-10pt} 
 
 \noindent\dotfill
 
  #3
  
   \noindent\dotfill 

\vspace{-10pt} 
 
 \noindent\hrulefill
 
 \vspace{-10pt} 
 
 \noindent\hrulefill
  \begin{quote} {\corcommstyle $>>$END$>>$ } \end{quote}
 \else
  #3
 \fi
}

\long\def\longcorrl #1<<#2||#3>>{
\if\visiblecomments y
  \begin{quote} {\corcommstyle $<<$COMMENT$>>$ {\color{red}#1\marginpar{!!}}\\$<<$OLD$<<$} \end{quote}
 
 {\color{red}

  #2
  
  }
  
  \begin{quote} {\corcommstyle ==NEW== } \end{quote}
  
    \noindent\hrulefill
 
\vspace{-10pt} 
 
 \noindent\hrulefill
 
 \vspace{-10pt} 
 
 \noindent\dotfill
 
  #3
  
   \noindent\dotfill 

\vspace{-10pt} 
 
 \noindent\hrulefill
 
 \vspace{-10pt} 
 
 \noindent\hrulefill
  \begin{quote} {\corcommstyle $>>$END$>>$ } \end{quote}
 \else
  #3
 \fi
}

\def\mlabel #1
{
  \if\visiblecomments y
     \marginpar[\flushright \bf \footnotesize #1]{\bf \footnotesize #1}
  \fi
}

\def\flabel #1
{
  \if\visiblecomments y
       \hbox{\bf\footnotesize #1}
  \fi
}

\def\corrq #1<<#2>>{
\if\visiblecomments y
  \begin{quote} {\corcommstyle $<<$COMMENT$>>$ {\color{red}#1}\marginpar{!!}\\$<<$BEG$<<$} \end{quote}
  \noindent\hrulefill
 
\vspace{-10pt} 
 
 \noindent\hrulefill
 
 \vspace{-10pt} 
 
 \noindent\dotfill

  #2
 
  \noindent\dotfill 

\vspace{-10pt} 
 
 \noindent\hrulefill
 
 \vspace{-10pt} 
 
 \noindent\hrulefill 
  \begin{quote} {\corcommstyle $>>$END$>>$ } \end{quote}
 \else
  #2
 \fi
}

\long\def\longcorrq #1<<#2>>{
\if\visiblecomments y
  \begin{quote} {\corcommstyle $<<$COMMENT$>>$ #1\marginpar{!!}\\$<<$BEG$<<$} \end{quote}
  \noindent\hrulefill
 
\vspace{-10pt} 
 
 \noindent\hrulefill
 
 \vspace{-10pt} 
 
 \noindent\dotfill

  #2

  \noindent\dotfill 

\vspace{-10pt} 
 
 \noindent\hrulefill
 
 \vspace{-10pt} 
 
 \noindent\hrulefill 
  \begin{quote} {\corcommstyle $>>$END$>>$ } \end{quote}
 \else
  #2
 \fi
}

\def\corrc #1<<>>{
\if\visiblecomments y
  \begin{quote} {\corcommstyle $<<$COMMENT$>>$ \color{red} #1\marginpar{!!}} \end{quote}
\fi
}

\def\corre #1<<#2||#3>>{
\if\visiblecomments y
  #3\marginpar{\corcommstyle #1}
 \else
  #3
 \fi
}

\long\def\longcorre #1<<#2||#3>>{
\if\visiblecomments y
  #3\marginpar{\corcommstyle #1}
 \else
  #3
 \fi
}

\def\corrs #1<<#2||#3>>{
\if\visiblecomments y
  #3\marginpar{\corcommstyle #2 $\rightarrow$ #3\\ #1}
 \else
  #3
 \fi
}

\def\corro #1<<#2||#3>>{
#2}

\def\corrn #1<<#2||#3>>{
#3}

\long\def\longcorro #1<<#2||#3>>{
#2}

\long\def\longcorrn #1<<#2||#3>>{
#3}

\long\def\underconstruction #1<<<#2>>>{
\if\visiblecomments y
  \begin{quote} {\corcommstyle $<<$UNDER CONSTRUCTION - BEGIN$>>$ #1\marginpar{!!}} \end{quote}
  #2
  \begin{quote} {\corcommstyle $>>$UNDER CONSTRUCTION - END$>>$ } \end{quote}
 \else
 \fi
}

\def\showcomments{
  \let\visiblecomments y
}

\def\hidecomments{
  \let\visiblecomments n
}

\showcomments
\title{Model rejection and parameter reduction via time series}

\author{Bree Cummins}
\address{Department of Mathematical Sciences \\
Montana State University\\
Bozeman, MT 59715}
\author{Tomas Gedeon}
\address{Department of Mathematical Sciences \\
Montana State University\\
Bozeman, MT 59715}
\author{Shaun Harker  }
\address{Department of Mathematics, Hill Center-Busch Campus\\
Rutgers, The State University of New Jersey \\
Piscataway, NJ  08854-8019, USA}
\author{Konstantin Mischaikow}
\address{Department of Mathematics, Hill Center-Busch Campus\\
Rutgers, The State University of New Jersey \\
Piscataway, NJ  08854-8019, USA}

\begin{document}

\maketitle

\begin{abstract}

We show how a graph algorithm for finding matching labeled paths in pairs of labeled directed graphs can be used to perform model validation for a class of dynamical systems including regulatory network models of relevance to systems biology.
In particular, we extract a partial order of events describing local minima and local maxima of observed quantities from experimental time-series data from which we produce a labeled directed graph we call the \emph{pattern graph} for which every path from root to leaf corresponds to a plausible sequence of events.
We then consider the regulatory network model, which can be itself rendered into a labeled directed graph we call the \emph{search graph} via techniques previously developed in computational dynamics.
Labels on the pattern graph correspond to experimentally observed events, while labels on the search graph correspond to mathematical facts about the model.
We give a theoretical guarantee that failing to find a match invalidates the model.
As an application we consider gene regulatory models for the yeast \textit{S. cerevisiae}.

\end{abstract}

\section{Introduction}
\label{sec:intro}
%!TEX root = ../PatternMatchPaper.tex

One of the fundamental challenges, as we move towards an era of data driven science, is how to make use of imprecise data to select or reject models and parameters that cannot be derived from first principles.
Motivated by problems from systems biology we address this challenge in the context of oscillatory data under the assumption that reasonable models  prescribe appropriate local behavior of trajectories.
We adopt the following strategy.
From experimental time series data we extract a partial order of events describing  minima and maxima of observed quantities.
On the modeling side, as a function of parameters, we construct a directed graph to catalogue the possible dynamics.
The main result of this paper is an efficient algorithm to identify if the model dynamics is capable of exhibiting sequences of minima and maxima that are consistent with the experimental data.
Failure can then be used  for model rejection or parameter reduction.

To provide more detail we consider a particular example.
High throughput  experimental technology is making the collection of time series of gene expression a routine process.
However, this data is noisy, often contains significant measurement error, is typically collected at a coarse time-scale, and generally is collected over a relatively short time span.
In an attempt to extract robust information from such data we focus on the ordering of extremal events.
This paper does not address the difficulty of detecting and eliminating spurious pairs of extrema in data -- a challenging problem in its own right.
Instead, we assume that a statistically valid procedure is used identify or impose time intervals during which a local maximum or minimum has occurred.
This renders the time series into a set of extremal events where the error bounds determine the  time intervals associated with the individual maxima and minima.
If two intervals do not overlap, then we can distinguish the relative timing between the associated events, but if they do overlap, we cannot.
Therefore, we represent relative timing as a partially ordered set (poset) that we call the \emph{poset of extrema} (see Definition~\ref{defn:PosetOfExtrema}).
We assume, however, that there is a linear temporal ordering along which the extrema occur and our lack of knowledge is  due to experimental constraints.
Consequently, we adopt the hypothesis that one of the linear extensions of the poset of extrema represents the correct sequence of events.

Because gene expression data is noisy and often collected at a coarse time scale, in practice there are many nodes in the poset of extrema that are not related.
Lack of relations lead to a multiplicative factors in the number of possible linear extensions.
Therefore, in general we expect that the set of all linear extensions will be large.
To organize the set of all linear extensions so that any particular linear extension can be accessed in a computationally efficient manner  we construct a labeled  directed acyclic graph that we call the \emph{pattern graph} (see Definition~\ref{defn:patterngraph}).
This directed graph is based on the ordering of the lattice of down sets of the poset of extrema.
The labeling acts on nodes and edges and is used to identify which variables are increasing, decreasing or have reached extrema.

The poset of extrema and the pattern graph represents the codification of the experimental data.
Viewed abstractly this is just a means of formally capturing potential temporal ordering of experimentally observable phenomena and therefore these ideas are potentially applicable to a wide variety of problems within and outside of the life sciences.

Returning to the example of gene regulation, we observe that this is an extremely complex multiscale process and thus it is not reasonable to  postulate a precise nonlinear model that describes its behavior.
However, as indicated above, we assume that we can make assumptions concerning the local qualitative behavior of the dynamics.
With this in mind we introduce the following notion.

\begin{defn}
\label{defn:dynamicalsystem}
A \emph{system of trajectories} $\mathcal{ST}$ on a space $X$ is a collection of continuous functions from closed intervals to $X$, i.e., $x \colon [a, b] \to X$ for some $a, b \in R$,  called \emph{trajectories}, such that: 
\begin{enumerate}
\item  the restriction of any trajectory to a smaller closed interval is again a trajectory; 
\item  a concatenation, i.e.\ continuous pasting, of trajectories is again a trajectory; 
\item  the time-translation, i.e.\ $x(t-\Delta t)$, of any trajectory is again a trajectory; and 
\item  every map $x \colon \{0\} \to X$ is a trajectory.
\end{enumerate}
\end{defn}

A heuristic description of how we employ this concept (see Section~\ref{section:patternmatch-searchgraph} for formal definitions) is as follows. 
We decomposed the phase space $X$ into a finite number of rectangular domains $\cX$ and consider a family of models for which trajectories within the domains are monotone and the trajectories on the boundaries between domains, called \emph{walls}, undergo at most one extremal event.
A system of trajectories that satisfies these conditions is called \emph{extrema-pattern-matchable} with respect to $\cX$.
The dynamics is then recorded as a labeled directed graph, called the \emph{search graph} (see Definition~\ref{defn:searchgraph}), where vertices correspond to domains and edges are determined by wall trajectories. 
The labeling acts on nodes and edges and is used to codify our knowledge as to which variables are increasing, decreasing, or have reached extrema.
Thus, the search graph formalizes the structure of the dynamics that can be expressed by a model that generates the system of trajectories.

The content  of this paper is the development of an efficient means of comparing the model dynamics against the experimentally observed dynamics.
The fundamental result  is Theorem~\ref{thm:AlignmentAlgorithm} that provides a polynomial time algorithm for matching maximal paths in the pattern graph, i.e.\ a specific ordering of extrema events that is compatible with the experimental data, to paths in the search graph, i.e.\ trajectories that are realizable by the model for the dynamics, where the matching preserves the labeling.
A simplistic description of the applicability of this result is as follows: if the algorithm fails to produce a matching, then this provides a guarantee that the dynamics incorporated in the system of trajectories is incapable of reproducing the sequences of extrema that are compatible with the data. 
As a consequence we can reject the associated model of dynamics.

Of course, to apply these ideas to realistic problems is much more challenging.
While we do not know a particular model that describes the observed dynamics we assume that the dynamics can be modeled by an unknown nonlinear system.
One of the fundamental lessons of the theory of dynamical systems is that structure of invariant sets of nonlinear systems can change dramatically as a function of parameters.
Thus, to achieve the claims of the title of this paper, we need  both a systematic method for generating parameterized models and their associated search graphs, and a robust finite characterization of dynamics that can be computed over all parameter values.
For this we make use of the recently developed Dynamics Signatures Generated by Regulatory Network (DSGRN) framework and software \cite{cummins2016combinatorial,DSGRNRepository}. 
%\corrc need Shaun's DOI for the code reference <<>>

The starting point for DSGRN is a regulatory network  $\textbf{RN}$ (see Definition~\ref{defn:rn}).
In the context of gene regulation this is an annotated directed graph in which the nodes represent genes, edges indicate the interaction between the genes, and the annotation indicates if the interaction involves activation or repression. 
Given a regulatory network with $N$ nodes and $|E|$ edges DSGRN represents dynamics occurring on the phase space $X=(0,\infty)^N$ where the dynamics is parameterized by an ${N+3\cdot |E|}$ dimensional set $Z \subset (0,\infty)^{N+3\cdot |E|}$.
In particular, DSGRN computes a finite decomposition of $Z$ where the individual regions are given by explicit semi-algebraic sets.
DSGRN represents parameter space via a undirected graph $\mathsf{PG}$, called the \emph{parameter graph}, where the nodes correspond to the above mentioned regions of $Z$ and edges provide adjacency information.
The dynamics is represented by a directed graph, called a \emph{state transition graph}, derived from a parameter dependent rectangular decomposition of $X$.
A fundamental fact is that as a function of parameters the state transition graph is constant over nodes of $\mathsf{PG}$, i.e., it does not change on the individual regions of the decomposition of parameter space.
The state transition graph can be large, and therefore DSGRN condenses the information into a directed acyclic graph $\mathsf{MG}$, called a \emph{Morse graph}.
The nodes of the Morse graph correspond to maximal recurrent subgraphs of the state transition graph and the edges indicate reachability, via the state transition graph, from one Morse node to another.
The output of DSGRN is called the DSGRN database, which is organized around the parameter graph $\mathsf{PG}$. 
In particular, for each node in $\mathsf{PG}$ the database provides the explicit semi-algebraic set in parameter space and the associated Morse graph $\mathsf{MG}$.

Observe that, as desired, DSGRN provides us with a finite  description of global dynamics over parameter space.
However, the dynamics is described in combinatorial terms and we would like to argue that we are comparing the ordering extrema of continuous trajectories against experimental data.
To make this comparison as transparent as possible we use the information encoded in the state transition graph to construct a particularly simple $\mathcal{ST}$ based on classical switching system models~\cite{DeJong2002,Edwards2000,Edwards2012,Edwards2014,Glass1972,Glass1973,Glass1978,Gouze2002}.

This paper is organized as follows. In Section~\ref{section:GraphTheory} we begin by recalling ideas from and establishing notation associated with graphs and posets. We  present the algorithms that underlie our approach to matching model dynamics with experimental data and provide worst case complexity bounds for these algorithms.

In Section~\ref{section:PatternMatching} we provide combinatorial formalizations of the experimental data, the dynamics of the models, and the relation that allows us to compare them. 
To be more specific, in Section~\ref{section:patterngraph} we show how experimental data can furnish a \emph{poset of extrema} $P$.
Interest in the set of all linear extensions leads us (by Theorem~\ref{thm:bijection}) to construct the down set graph of the poset of extrema.
We label each down set according to whether a function that has experienced the events in the down set but not the events not included in the down set is increasing or decreasing in each variable.
We label the edges in the down set graph (which are of the form $A \to A \cup \{p\}$) according to the extremal event associated with $p \in P$.
We also introduce self-edges that are labeled as not experiencing any extremal events.
We call the resulting labeled directed graph the \emph{pattern graph}.

In Section~\ref{section:patternmatch-searchgraph}, we describe a class of dynamical models for which we can characterize possible trajectories in a combinatorial manner via a \emph{domain graph}.
A domain graph discretizes a dynamical system by giving a finite set of domains separated by codimension-1 walls.
Vertices in the domain graph correspond to domains, and edges correspond to flow from one domain to another via a wall.
We label the vertices of the domain graph according to whether the coordinate functions $x_i(t)$ are increasing, decreasing, or possibly both, in the associated domains.
We label the edges of the graph according to which local minima or local maxima could occur on the associated walls.
We call the resulting labeled directed graph a \emph{search graph}.

In Section~\ref{section:matchingtheorem}, we present a matching relation between the pattern graph and the search graph.
We prove (Theorem~\ref{thm:NoFalseNegatives}) that if there does not exist a match between a path from root to leaf of the pattern graph and a path in the search graph, then the dynamical model underlying the search graph is incompatible with the experimental observations leading to the pattern graph.
Theorem~\ref{thm:AlignmentAlgorithm} shows that we can decide whether or not such a match exists in polynomial time.

Finally, in Section~\ref{section:Application} we show how these ideas can be applied.
We begin in Section~\ref{section:DSGRN} with a brief review of the mathematical structure underlying DSGRN.
In Section~\ref{section:applicationpatterngraph} we provide a simple example of how one can pass from experimental time series data to a labeled pattern graph. 
Finally, in Section~\ref{section:ApplicationResults} we apply these techniques to a simple wavepool model~\cite{orlando2008global} for  the metabolic cycle in \textit{S. cerevisiae}. 
Courtesy of the Haase lab~\cite{leman2014analyzing} we have experimental time-series data for  mRNA sequences associated with the genes SWI4, HCM1, NDD1, and YOX1 collected at time intervals of 5 minutes (see Figure~\ref{fig:dataset}).
We take a biologically implausible model and show that our proposed techniques reject it and we take a model that is biologically acceptable and use our techniques to greatly constrain relations between parameters.

\section{Graph Theory and Algorithms}
\label{section:GraphTheory}
%!TEX root = ../PatternMatchPaper.tex

\subsection{Matching Paths in Labeled Graphs}
\label{subsection:MatchingPaths}

\begin{defn}
\label{defn:Kleene}
Given a finite set $\Sigma$, we denote by $\Sigma^n$ the set of $n$-tuples consisting of elements of $\Sigma$.
We denote by $\Sigma^*$ the set of all finite tuples of elements of $\Sigma$, i.e. $\Sigma^* = \bigcup_{n=0}^\infty \Sigma^n$.
\end{defn}

For the purpose of this paper  a \emph{directed graph} $G = (V,E)$ consists of a finite set of \emph{vertices} $V$ and \emph{edges} $E \subset V\times V$.
A \emph{path} in $G$ from $s\in V$ to $t \in V$ is a finite sequence of vertices $(s = v_1, v_2, \cdots, v_n = t)$ such that $v_i \in V$ and $(v_i, v_{i+1}) \in E$.
We denote the set of all such paths by $G[s \leadsto t]$.

\begin{defn}
\label{defn:LabeledDirectedGraph}
A \emph{labeled directed graph} $G$ is a quadruple $(V, E, \Sigma, \ell)$ where $V$ and $E$ denote the vertices and edges of $G$,  $\Sigma$ is a finite set called \emph{labels}, and $\ell : V \cup E \to \Sigma$ is called a \emph{labeling function}.
Given a path $p = (v_1,\ldots, v_n)$ in $G$, the associated \emph{labeling} is defined to be 
\[ 
L(p) := (\ell(v_1), \ell((v_1,v_2)), \ell(v_2), \cdots, \ell((v_{n-1},v_n)), \ell(v_n)) \in \Sigma^*.
\]
\end{defn}

\begin{defn}
\label{defn:MatchingPaths}
A \emph{matching relation} between two labeled directed graphs $G = (V, E, \Sigma, \ell)$ and $G' = (V', E', \Sigma', \ell')$ is a relation $\sim$ between the label sets $\Sigma$ and $\Sigma'$.
The labels $a \in \Sigma$ and $b \in \Sigma'$ \emph{match} if $a \sim b$.
We extend the matching relation $\sim$ onto the tuples of labels $\Sigma^*$ and $\Sigma'^*$ via \[ (a_1, a_2, \cdots, a_n) \sim (b_1, b_2, \cdots, b_m) \mbox{ iff } n=m \mbox{ and for $1 \leq i \leq n$, } a_i \sim b_i. \]
Given a matching relation, a path $p=(v_1,\cdots,v_n)$ in $G[s \leadsto t]$, and a path $p'=(v'_1,\cdots,v'_m)$ in $G'[s' \leadsto t']$, we say that \emph{$p$ matches $p'$} and write $p \sim p'$ whenever $L(p) \sim L(p')$.
Note that we are using the same symbol $\sim$ to refer to three matching relations: between $\Sigma$ and $\Sigma'$, between $\Sigma^*$ and $\Sigma'^*$, and between paths in $G[s \leadsto t]$ and paths in $G'[s' \leadsto t']$.
\end{defn}

\begin{defn}
Let $\sim$ be a matching relation between two labeled directed graphs $G = (V, E, \Sigma, \ell)$ and $G' = (V', E', \Sigma', \ell')$.
Suppose $s,t \in V$ and $s',t' \in V'$.
The \emph{alignment problem} $\mathsf{Alignment}(G, G', \sim, (s,t), (s', t'))$ is the decision problem of determining if there is a pair of paths $p \in G[s \leadsto t]$ and $p'\in G'[s' \leadsto t']$ such that $p \sim p'$.
\end{defn}

\begin{thm}
\label{thm:AlignmentAlgorithm}
There exist polynomial time algorithms for the following decision problems:
\begin{enumerate}
\item Let $s,t \in V$, $s',t' \in V'$. Decide $\mathsf{Alignment}(G, G', \sim, (s,t), (s', t'))$
\item Let $s,t \in V$. Decide $\exists s',t' \in V'\; \mathsf{Alignment}(G, G', \sim, (s,t), (s', t'))$
\item Let $s,t \in V$. Decide $\exists s' \in V'\; \mathsf{Alignment}(G, G', \sim, (s,t), (s', s'))$
\end{enumerate}
\end{thm}

We postpone the proof of Theorem~\ref{thm:AlignmentAlgorithm} to Section~\ref{section:Algorithms}, where we give explicit algorithms.

%!TEX root = ../PatternMatchPaper.tex

\subsection{Down Set Graph of a Poset}

\begin{defn}
A \textit{poset} $(P,<)$ is a set $P$ equipped with a transitive, irreflexive relation $<$ called a \emph{partial order}.
A \textit{linear extension} of $<$ is a total order $<'$ which extends $<$, i.e. for all $p_0,p_1\in P$, $p_0<p_1$ implies $p_0<'p_1$.
\end{defn}

\begin{defn}
Let $(P, <)$ be a poset.
A \emph{down set} of $P$ is a subset $A\subset P$ such that for all $p, q \in P$, $p < q$ and $q \in A$ implies $p\in A$.
The collection of down sets of $P$ is denoted by $O(P)$.
\end{defn}

\begin{defn}
\label{defn:downsetgraph}
Let $(P, <)$ be a finite poset.
The \textit{down set graph} of $(P,<)$, denoted $P_D$, is the directed graph $(O(P),F)$ with vertices $O(P)$ and edges $A \rightarrow A'$ if and only if $A\subsetneq A'$ and there does not exist $A''\in O(P)$ such that $A\subsetneq A''\subsetneq A'$.
\end{defn}

\begin{rem}
\label{rem:downsetgraph}
If $A,A'\in O(P)$ and there exists an edge $A \rightarrow A'$ in $P_D$, then $A' = A\cup\setof{p}$ where $p\not\in A$ and $q<p$ implies that $q\in A$.
\end{rem}

For completeness, we include an algorithm for constructing the down set graph of a poset in Section~\ref{section:constructionofdownsetgraph}.

\begin{thm}
\label{thm:bijection}
Given a finite poset $(P,<)$, the associated down set graph $P_D = (O(P),F)$ is a directed acyclic graph with a unique root $\emptyset$ and a unique leaf $P$.
Moreover, there is a bijection between the paths in $P_D$ from the root $\emptyset$ to the leaf $P$ and the linear extensions of $<$.
\end{thm}

\begin{proof}
The directed acyclic property is inherited from the definition via proper set inclusion.
$\emptyset$ and $P$ are the unique root and leaf since $\emptyset$ and $P$ are the unique maximal and minimal elements in $O(P)$, respectively.
Now we show the moreover part.
Let $<'$ be a linear extension of $<$.
Suppose $P = \{p_1, p_2, \cdots, p_n\}$ where the indexing has been chosen so that $p_1 <' p_2 <' \cdots <' p_n$.
Define $P_k := \{p_1, p_2, \cdots, p_k\}$.
Then $\emptyset \to P_1 \to P_2 \to \cdots \to P_{n-1} \to P_{n} = P$ gives a path in $P_D$ from root to leaf unique to $<'$.
Now the converse.
Suppose $\emptyset \to P_1 \to P_2 \to \cdots \to P_{n-1} \to P_{n} = P$ is a path from root to leaf in $P_D$.
We claim that $P_k \setminus P_{k-1}$ must be a singleton for each $k$.
Suppose otherwise.
Then let $a, b \in P_k \setminus P_{k-1}$ such that $a \neq b$.
Without loss, assume either $a < b$ or $a$ and $b$ are incomparable.
Then $P_{k-1}\cup\{a\}$ is a down set and $P_{k-1} \subsetneq P_{k-1}\cup\{a\} \subsetneq P_k$ which by Definition~\ref{defn:downsetgraph} contradicts $P_{k-1} \to P_{k}$.
Accordingly, let $p_k = P_{k} \setminus P_{k-1}$ for $k=1,\cdots,n$, and see that this sequence of elements completely characterizes the path from root to leaf in $P_D$.
Define the total order $<'$ via $p_1 <' p_2 <' \cdots <' p_n$; since $p_k < p_{k+1}$ holds for all $k$, $<'$ is a linear extension of $<$.
Thus a path from root to leaf in $P_D$ uniquely determines a linear extension $<'$ of $<$.
\end{proof}

%!TEX root = ../PatternMatchPaper.tex

\subsection{Algorithms} \label{section:Algorithms}
\subsubsection{Alignment Problem}
\label{subsection:AlignmentAlgorithm}

\begin{defn}\label{defn:alignmentgraph}
Let $\sim$ be a matching relation between two labeled directed graphs $G = (V, E, \Sigma, \ell)$ and $G' = (V', E', \Sigma', \ell')$.
The \emph{alignment graph} $\mathsf{AlignmentGraph}(G, G', \sim)$ is defined to be the directed graph $(V'', E'')$ given by
\[ V'' = \{ (v,v') \in V \times V' : \ell(v) = \ell'(v') \}, \]
\[ E'' = \{ (e, e') \in E \times E' : \ell(e) = \ell'(e') \}.\]
\end{defn}

The alignment graph $\mathsf{AlignmentGraph}(G, G', \sim)$ is a subset of the product graph $G\times G'$, and hence it has at most $\vert V \vert \vert V' \vert$ vertices, $\vert E \vert \vert E' \vert$ edges.

The following proposition follows immediately from the construction of the alignment graph and the definition of matching paths:

\begin{prop}
Paths in $\mathsf{AlignmentGraph}(G, G', \sim)$ are in one-to-one correspondence with pairs of matching paths in $G$ and $G'$.
In particular, $p''=((v_1, v'_1), (v_2, v'_2), \cdots, (v_n, v'_n))$ is a path in the alignment graph if and only if $p = (v_1, v_2, \cdots, v_n)$ and $p' = (v'_1, v'_2, \cdots, v'_n)$ are a pair of matching paths in $G$ and $G'$ respectively.
\end{prop}

It immediately follows that the alignment problem is equivalent to a reachability query in the alignment graph:

\begin{prop}
\label{prop:AlignmentGraphProperty}
The following are equivalent:
\begin{enumerate}
\item $\mathsf{Alignment}(G, G', \sim, (s,t), (s', t'))$
\item $\mathsf{AlignmentGraph}(G, G', \sim)\left[(s,s') \leadsto (t, t')\right] \neq \emptyset$
\end{enumerate}
\end{prop}

\begin{prop}
\label{prop:AlignmentGraphConstruction}
If the cost of checking whether labels match is constant, then $\mathsf{AlignmentGraph}(G, G', \sim)$ can be constructed in $O(\vert V \vert \vert V' \vert + \vert E \vert \vert E' \vert)$ time.
\end{prop}

\begin{proof}
The vertices of the alignment graph may be determined by checking for each element of $(v, v') \in V \times V'$ whether $\ell(v) = \ell(v')$.
The edges of the alignment graph may be determined by checking for each $(e, e') \in E \times E'$ whether $\ell(e) = \ell(e')$.
The result follows.
\end{proof}

\begin{prop}
\label{prop:ReachabilityAlgorithm}
Let $G = (V, E)$ be a directed graph and let $s\in G$.
Then, the set $\mathsf{Reachable}(G,s) := \{ t \in V : G[s \leadsto t] \neq \emptyset \}$ can be computed in $O(\vert V \vert + \vert E \vert)$ time.
\end{prop}

\begin{proof}
Depth or breadth first search of $G$ beginning at $s$ will find all vertices reachable from $s$ in time linear in the number of vertices and edges of $G$ \cite{Cormen:2001:IA:580470}.
For completeness we provide a standard depth-first-search algorithm as $\mathsf{Reachable}(G,s)$ in Algorithm~\ref{alg:alignment}.
\end{proof}

\noindent\emph{Proof of Theorem~\ref{thm:AlignmentAlgorithm}.}
We show that the procedures $\mathsf{Match}$, $\mathsf{PathMatch}$, and $\mathsf{CycleMatch}$ of Algorithm~\ref{alg:alignment} are polynomial time algorithms which decide the decision problems (1), (2), and (3), respectively.
Let $G'' = \mathsf{AlignmentGraph}(G, G', \sim)$.
By Proposition~\ref{prop:AlignmentGraphProperty}, we may rewrite (1), (2), (3) as
\begin{enumerate}
\item Let $s,t \in V$, $s',t' \in V'$. Decide if $(t,t') \in \mathsf{Reachable}(G'', (s,s'))$
\item Let $s,t \in V$. Decide $\exists s',t' \in V'\; (t,t') \in \mathsf{Reachable}(G'', (s,s'))$
\item Let $s,t \in V$. Decide $\exists s' \in V'\; (t,s') \in \mathsf{Reachable}(G'', (s,s'))$
\end{enumerate}

The correctness of $\mathsf{Match}$ for deciding (1) is now immediate.
For (2) and (3), we recognize we can handle the outermost $\exists s' \in V'$ algorithmically via a \textbf{for} loop over $s' \in V'$.
The algorithms $\mathsf{PathMatch}$, and $\mathsf{CycleMatch}$ result.
This gives correctness.

To see that the algorithms are polynomial time, we refer to Proposition~\ref{prop:AlignmentGraphConstruction} and Proposition~\ref{prop:ReachabilityAlgorithm}.
In particular, given these it is straightforward to verify (defining $\vert G \vert = \vert V \vert + \vert E \vert$,) that $\mathsf{Reachable}$ executes in worst-case $O(\vert G \vert)$ time, $\mathsf{Match}$ executes in worst-case $O(\vert G \vert \vert G' \vert)$ time, and $\mathsf{PathMatch}$ and $\mathsf{CycleMatch}$ execute in worst-case $O(\vert G \vert \vert G' \vert \vert V' \vert)$ time.
\qedsymbol\\

See Figure~\ref{fig:matching_path} for an example of a $\mathsf{PathMatch}$ along with the corresponding path giving reachability in the alignment graph.

\begin{algorithm}[h]
\caption{Alignment Problem}
\label{alg:alignment}

\begin{multicols}{2}

\begin{algorithmic}
\Procedure{Reachable}{$G$, $s$}
\State Push $s$ onto stack $S$.
\While{$S$ is not empty}
  \State Pop $u$ from stack $S$.
  \State $R \gets R \cup \{u\}$
  \State $A \gets \{ v : (u, v) \in E \}$
  \For{$v \in A$}
    \If{ $v \notin R$ }
      \State Push $v$ into stack $S$
    \EndIf
  \EndFor
\EndWhile
\State \textbf{return} $R$
\EndProcedure
\end{algorithmic}

\vspace{8pt}

\begin{algorithmic}
\Procedure{Match}{$G, G', (s, t), (s', t')$}
\State $G'' \gets \mathsf{AlignmentGraph}(G, G', \sim)$
\State $R \gets \mathsf{Reachable}(G'', (s,s'))$
\If{ $(t,t') \in R$ }
  \State \textbf{return} True
\Else
  \State \textbf{return} False
\EndIf
\EndProcedure
\end{algorithmic}

\columnbreak

\begin{algorithmic}
\Procedure{PathMatch}{$G$, $G'$, $(s, t)$}
\State $G'' \gets \mathsf{AlignmentGraph}(G, G', \sim)$
\For{$s' \in V'$}
  \State $R \gets \mathsf{Reachable}(G'', (s,s'))$
  \For{$(v,v') \in R$}
    \If{ $ v = t$ }
      \State \textbf{return} True
    \EndIf
  \EndFor
\EndFor
\State \textbf{return} False
\EndProcedure
\end{algorithmic}

\vspace{8pt}

\begin{algorithmic}
\Procedure{CycleMatch}{$G$, $G'$, $(s, t)$}
\State $G'' \gets \mathsf{AlignmentGraph}(G, G', \sim)$
\For{$s' \in V'$}
  \State $R \gets \mathsf{Reachable}(G'', (s,s'))$
  \For{$(v,v') \in R$}
    \If{ $ v = t$ and $v' = s'$ }
      \State \textbf{return} True
    \EndIf
  \EndFor
\EndFor
\State \textbf{return} False
\EndProcedure
\end{algorithmic}
\end{multicols}

\end{algorithm}

\subsubsection{Construction of Down Set Graph}
\label{section:constructionofdownsetgraph}

\begin{algorithm}[h]
\caption{Down Set Graph}
\label{alg:posettodownset}
\begin{algorithmic}
\Procedure{PosetToDownSetGraph}{$P$}
\State Let $S$ be an empty stack
\State $V \gets \emptyset$
\State $E \gets \emptyset$
\State Push $\mathsf{MaximalElementsOf}(P)$ onto $S$
\While{$S$ is not empty}
  \State Pop $I$ from $S$
  \State $V \gets V \cup \{I\}$
  \For{$v \in I$}
    \State $I' \gets \mathsf{MaximalElementsOf}((I \cup \mathsf{Predecessors}(v))\setminus \{v\})$
    \If{ $I' \notin V$}
      \State Push $I'$ onto $S$
    \EndIf
    \State $E \gets E \cup \{(I', I)\}$
  \EndFor
\EndWhile
\State \textbf{return} $(V,E)$
\EndProcedure
\end{algorithmic}
\end{algorithm}

Recall that given a poset $P$ a subset $I\subset P$ is \emph{independent} if no two elements of $I$ are comparable.

\begin{prop}
Algorithm~\ref{alg:posettodownset} computes the down set graph of a poset $P$. 
\end{prop}

\begin{proof}
We first note there is a one-to-one correspondence between down sets of a poset and the independent sets of a poset.
In particular given a down set $D$ we can associate an independent set $I = \mathsf{MaximalElementsOf}(D)$, and given an independent set $I$ we can associate a down set $D = \mathsf{Downset}(I) := \{ p \in P : p \leq q \mbox{ for some } q \in I\}$.
It is straightforward to see that this is one to one.
Now consider the following recursively defined function:

\[ f(D) := \{ (D, D \setminus \{v\} ) : \mbox{for $v \in \mathsf{MaximalElementsOf}(D)$} \} \cup \bigcup_{v \in \mathsf{MaximalElementsOf}(D)} f(D \setminus \{v\}) \]

See first that the recursion terminates since each recursive function call operates on a smaller set.
Notice that if the function operates on a down set then removing the maximal vertices again results in down sets -- in fact, precisely the adjacent down sets in the down set graph.
Hence $E = f(P)$ is the set of edges in the down set graph.
Writing this recursion in terms of independent sets, we have
\begin{align*}
g(I)  :=  & \setof{ (I, \mathsf{MaximalElementsOf}(\mathsf{Downset}(I) \setminus \{v\}) ) : \text{for $v \in I$} } \\
& \quad \cup \bigcup_{v \in I} g(\mathsf{MaximalElementsOf}(\mathsf{Downset}(I) \setminus \{v\}))
\end{align*}
Now from 
\[
\mathsf{MaximalElementsOf}(\mathsf{Downset}(I)\setminus \{v\}) = \mathsf{MaximalElementsOf}((I \cup \mathsf{Predecessors}(v))\setminus \{v\})
\] 
the correctness of Algorithm~\ref{alg:posettodownset} follows: it is just an implementation of this recursion which prevents some redundant recursion paths to save time (by storing them in $V$).
\end{proof}

Algorithm~\ref{alg:posettodownset} does not run in polynomial time in general, yet it does for the special case of interest for the application of this paper. In particular, as we will describe in the next section, we will consider posets for which each element is associated to one of a small number of variables $x_1, x_2, \cdots, x_d$, and all elements in the poset associated to the same variable are comparable. Moreover, as we discuss in Section~\ref{section:Application}, associated to poset elements are time intervals determining the partial order such that one time interval $(a,b)$ compares less than another time interval $(c, d)$ iff $b \leq c$. Under these assumptions, the \emph{incomparability graph of the poset $P$} (i.e. the graph with vertices $P$ and edges $u \leftrightarrow v$ whenever $u$ and $v$ are incomparable) is an \emph{interval graph} \cite{fulkerson1965incidence}, which is a special kind of \emph{chordal graph}. A chordal graphs with $n$ vertices has at most $n$ maximal cliques \cite{rose1976algorithmic, gavril1972algorithms}. From these considerations we get the following bound:

\begin{prop}
\label{prop:downsetcomplexity}
Assume that $P$ is a finite poset. Let $d$ be the cardinality of the maximum independent set in $P$. Assume that the incomparability graph of $P$ is chordal. Then the down set graph $D_P$ has at most $2^d n$ vertices, and for fixed $d$, Algorithm~\ref{alg:posettodownset} executes in polynomial time.
\end{prop}

\section{Matching Posets of Extrema against Computational Dynamics Models}
\label{section:PatternMatching}
% \input{texfiles/3-dynamics-overview}
%!TEX root = ../PatternMatchPaper.tex

\subsection{Labeled Directed Graphs from Posets of Extrema}
\label{section:patterngraph}

Assume that we can measure $N$ variables over a time interval $[0,T]$ for the system that we are interested in modeling.
If the quantities of these variables change continuously, then there exists a continuous function $x : [0,T] \to \R^N$ that represents the dynamics.
We will assume that over this time interval each variable attains finitely many local extrema.
As discussed in the introduction, in applications we can only sample the system at finite time intervals and the measurements will be subject to noise.
We use the following structure to codify the possible orderings of maxima and minima of the coordinates $x_i$ of $x$.
\begin{defn}
\label{defn:PosetOfExtrema}
A \emph{poset of extrema} $(P,<_\tau;\mu)$ is a finite poset $(P, <_\tau)$ equipped with a surjective function $\mu : P \to \{ -, m, M \}^N$ that satisfies the following conditions.
For $n=1,\ldots, N$, define $P_n = \{ p : \mu(p)_n \in \{m, M\} \}$.  
\begin{enumerate}
\item $P = \bigcup_{n=1}^N P_i$.
\item  If $n\neq j$, then $P_n\cap P_j = \emptyset$.
\item For each $n$, $P_n \subset P$ is totally ordered by $<_\tau$.
\item Let $u,v\in P_n$. If $u<_\tau v$ and $\mu(u) = \mu(v)$, then there exists $w \in P_n$ such that $u <_\tau w <_\tau v$ and $\mu(u)\neq \mu(w)$. 
\end{enumerate}
\end{defn}

It is worth commenting on the rationale behind Definition~\ref{defn:PosetOfExtrema}.
The poset of extrema is designed to capture orderings with respect to time of minima and maxima of $d$ variables.
The symbols $-,m,$ and $M$ stand for not an extremum, local minimum, and local maximum, respectively, and in applications the ordering $<_\tau $ respects the direction of time.
Condition (1) implies that every vertex of $P$ is associated with an extremal event. 
Condition (2) implies that each vertex is associated to an extremal event of precisely one variable, i.e., $\mu(p)_n = -$ for all but precisely one $n \in \{1,2,\cdots,N\}$.
The assumption that each $P_n$ is totally ordered with respect to $<_\tau$ implies that for each variable the ordering (with respect to time) of the minima and maxima is known; any ambiguity arises from comparing across variables.
The final condition prevents a variable experiencing two local maxima or two local minima consecutively.
Note that this is an assumption about the sampling frequency of the experiment.

Returning to the unknown function $x$ that represents the dynamics, one expects, generically, that  the maxima and minima of the coordinates $x_n$ occur at different times. 
In the context of the poset of extrema $(P,<_\tau)$ we interpret this to mean that the true dynamics corresponds to linear extension of $<_\tau$.
Since, given the data, the linear extension is unknown we consider any linear extension to be a plausible sequence of events.
Our goal is to use the the machinery of Section~\ref{subsection:MatchingPaths} in order to search for linear extensions of $P$ and thus we construct, following Theorem~\ref{thm:bijection}, the down set graph $P_D$ of $P$, which exhibits a one-to-one correspondence between paths from root to leaf and linear extensions of $P$.

In order to produce a labeled directed graph suitable for pattern matching algorithms, we require labels on the vertices and edges  of $P_D$. 
We make use of a particular set of labels 
\[
\Sigma_{\text{ext}} :=  \setof{I, D, *, -,m,M}
\]
called the \emph{extrema labels} which are intended to carry the following information:
\begin{itemize}
\item $I$: increasing,
\item $D$: decreasing,
\item $m$: minimum,
\item $M$: maximum,
\item $-$: transitioning,
\item $*$: lack of knowledge.
\end{itemize}

 \begin{defn}
\label{defn:patterngraph}
Let $(P,<_\tau ;\mu)$ be a poset of extrema with down set graph $P_D = (O(P), E)$.
The \emph{pattern graph} $\cP$ induced by the poset of extrema $P$ is the labeled directed graph $(O(P), E \cup \{(A,
A): A \in O(P)\}, \Sigma_{\text{ext}}^N, \ell)$, where the labeling of the vertices is given by
\[ \ell(A)_n = 
\begin{cases}
 I \mbox{ if $\ell(\max (P_n \cap A))_n = m$ or $\ell(\min (P_n \setminus A))_n = M$} \\
D \mbox{ if $\ell(\max (P_n \cap A))_n = M$ or $\ell(\min (P_n \setminus A))_n = m$} \\
 * \mbox{ otherwise,}\end{cases}\]
and the labeling of the edges is defined by
\[ \ell(A \to A) := (-,-,\cdots,-)\]
and (see Remark~\ref{rem:downsetgraph})
\[ \ell(A \to A \cup \setof{p}) := \mu(p),\quad p\in P.\]
 \end{defn}

Although the pattern graph is (trivially) cyclic due to the presence of self-edges, we will continue to refer to the root and leaf nodes of $\cP$ as $\emptyset$ and $P$ respectively.

\begin{figure}
   \begin{center}
    \begin{tikzpicture}[main node/.style={ellipse,fill=white!20,draw},scale=1.5]
      \node[main node] (Sm) at (0,0) {$x_1$ min};
      \node[main node] (SM) at (0,-1) {$x_1$ max};
      \node[main node] (Ym) at (2.3,0) {$x_2$ min};
      \node[main node] (YM) at (2.3,-1) {$x_2$ max};

      \path[->,>=angle 90,thick]
      (Sm) edge[shorten <= 3pt, shorten >= 3pt] node[] {} (SM)
      (Sm) edge[shorten <= 3pt, shorten >= 3pt] node[] {} (YM)
      (Ym) edge[shorten <= 3pt, shorten >= 3pt] node[] {} (SM)
      (Ym) edge[shorten <= 3pt, shorten >= 3pt] node[] {} (YM)
      ;

      \node[main node] (DD1) at (6,2) {0: DD};
      \node[main node] (ID1) at (5,1) {1: ID};
      \node[main node] (DI1) at (7,1) {2: DI};
      \node[main node] (II) at (6,0) {3: II};
      \node[main node] (ID2) at (7,-1) {4: ID};
      \node[main node] (DI2) at (5,-1) {5: DI};
      \node[main node] (DD2) at (6,-2) {6: DD};

      \path[->,>=angle 90,thick]
      (DD1) edge[shorten <= 3pt, shorten >= 3pt] node[near end,above=5pt] {m-} (ID1)
      (DD1) edge[shorten <= 3pt, shorten >= 3pt] node[near end,above=5pt] {-m} (DI1)
      (DD1) edge[loop left, shorten <= 3pt, shorten >= 3pt] node[] {-~-} (DD1)
      (ID1) edge[shorten <= 3pt, shorten >= 3pt] node[near start,below=5pt] {-m} (II)
      (DI1) edge[shorten <= 3pt, shorten >= 3pt] node[near start,below=5pt] {m-} (II)
      (ID1) edge[loop left, shorten <= 3pt, shorten >= 3pt] node[] {-~-} (ID1)
      (DI1) edge[loop right, shorten <= 3pt, shorten >= 3pt] node[] {-~-} (DI1)
      (II) edge[shorten <= 3pt, shorten >= 3pt] node[near end,above=5pt] {M-} (DI2)
      (II) edge[shorten <= 3pt, shorten >= 3pt] node[near end,above=5pt] {-M} (ID2)
      (II) edge[loop left, shorten <= 3pt, shorten >= 3pt] node[] {-~-} (II)
      (ID2) edge[shorten <= 3pt, shorten >= 3pt] node[near start,below=5pt] {M-} (DD2)
      (DI2) edge[shorten <= 3pt, shorten >= 3pt] node[near start,below=5pt] {-M} (DD2)
      (ID2) edge[loop right, shorten <= 3pt, shorten >= 3pt] node[] {-~-} (ID2)
      (DI2) edge[loop left, shorten <= 3pt, shorten >= 3pt] node[] {-~-} (DI2)
      (DD2) edge[loop left, shorten <= 3pt, shorten >= 3pt] node[] {-~-} (DD2)
      ;

      \end{tikzpicture}
    \end{center}\caption{Left: Example poset of extrema with four extrema. Right: Associated pattern graph.}\label{fig:earlypoexample}
 \end{figure}
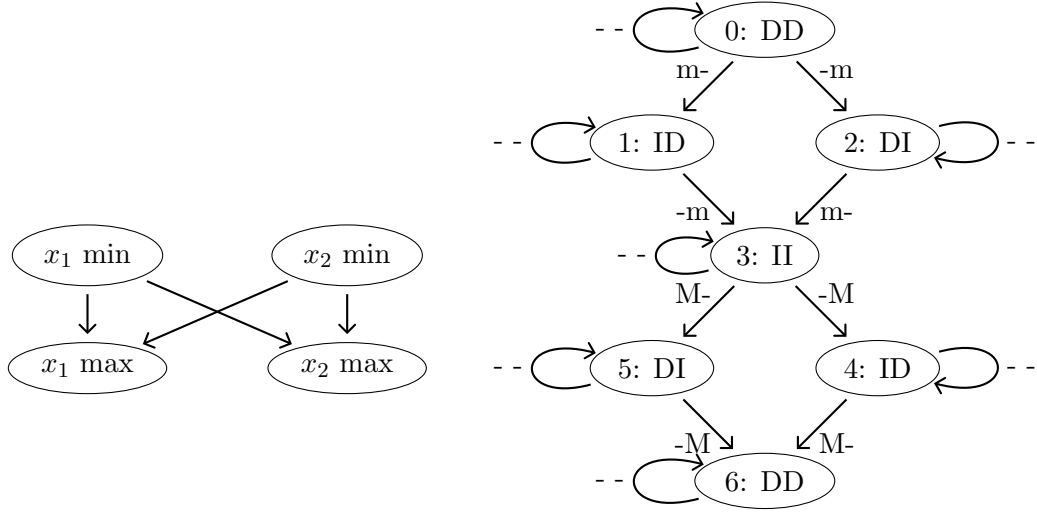

We give an example of a poset of extrema and the associated pattern graph using two variables $x_1$ and $x_2$, i.e. $N=2$.  Later on, we shall relate this example to a yeast dataset~\cite{leman2014analyzing} that is discussed in Section~\ref{section:applicationpatterngraph}. For now, assume that $x_1$ and $x_2$ first attain minima, then later attain maxima, but that  the timing of the  minima cannot be distinguished, and neither can the maxima. This leads to the poset in Figure~\ref{fig:earlypoexample} (left). The associated pattern graph is in Figure~\ref{fig:earlypoexample} (right). The down sets of the poset of extrema are mapped to integers via 
\begin{align*} &0 \leftrightarrow \emptyset; \; 1 \leftrightarrow \{x_1 \min\}; \; 2 \leftrightarrow \{x_2 \min\}; \; 3 \leftrightarrow \{x_1 \min, x_2 \min\}; \; 
\\ &4 \leftrightarrow \{x_1 \min, x_2 \min, x_2 \max\}; \; 5 \leftrightarrow \{x_1 \min, x_2 \min, x_1 \max\}; \; 
\\ &6 \leftrightarrow \{x_1 \min, x_2 \min, x_1 \max, x_2 \max\}. \end{align*}

%!TEX root = ../PatternMatchPaper.tex

\subsection{Labeled Directed Graphs From Computational Dynamics}
\label{section:patternmatch-searchgraph}

In this section we develop the notion of a search graph, a labeled directed graph suitable for pattern matching sequence of extrema for models arising in computational dynamics.

\begin{defn}
\label{defn:cubicaldecomposition}
Let $X = (0,\infty)^N$.
Suppose for each $n \in \{1,\cdots,N\}$ we have a finite set $\Theta_n \subset (0, \infty)$.
The 
\emph{rectangular decomposition $\cX$ of $X$ induced by the thresholds $(\Theta_1, \Theta_2, \cdots, \Theta_N)$} 
is the partition of $X$ into \emph{cells} $\cX$ corresponding to the classification of each point in $X$ according to whether or not it is contained in, above, or below each of the hyperplanes in $\{ \{ x : x_n = \theta \} : n \in \{1,\cdots, N\} \mbox{ and } \theta \in \Theta_n\}$.
It follows that for each cell $\sigma \in \cX$ is a product of intervals
\[ \sigma = \Pi_{n=1}^N I_n, \] where for each $n$,  $I_n \in \{ (0, \theta_1), (\theta_1,\theta_2), (\theta_2,\infty),  [\theta_1,\theta_1] \}$ for some $\theta_1, \theta_2 \in \Theta_n$.
Accordingly, each cell is homeomorphic to an open ball of some dimension, which we call the dimension of the cell.
We denote $k$-dimensional cells $\cX_k$.
We call the cells in $\cX_N$ \emph{domains} and we call the cells in $\cX_{N-1}$ \emph{walls}.
Two domains are said to be \emph{adjacent} if the intersection of their closures contains a wall.
We denote by $X_k$ the union of all $k$-dimensional cells, i.e. $X_k := \bigcup_{\sigma \in \cX_k} \sigma$.
\end{defn}

\begin{defn}
\label{defn:domaingraph}
Consider $X = (0,\infty)^N$ with rectangular decomposition $\cX$ and a system of trajectories $\mathcal{ST}$  on $X$.
A trajectory $x:[t_0,t_1] \to X$ is a \emph{domain trajectory} if $x([t_0,t_1]) \subset \xi$ for some domain $\xi \in \cX_N$.
A trajectory $x : [t_0, t_1] \to X$ is a \emph{wall trajectory} from a domain $\xi$ to a domain $\xi'$ if there exists a wall $\sigma \in \cX_{N-1}$ such that $x([t_0,t_1]) \subset \xi \cup \sigma \cup \xi'$, and $x^{-1}(\xi) < x^{-1}(\sigma) < x^{-1}(\xi')$ (in the sense of comparing sets, i.e. $A < B$ iff $\forall a\in A,b\in B, a<b$).
The \emph{domain graph generated by $\mathcal{ST}$ on $\cX$} is the directed graph where the vertices are domains and there is an edge $\xi \to \xi'$ if and only if there exists a wall trajectory from $\xi$ to $\xi'$.
\end{defn}

Definition~\ref{defn:domaingraph} indicates how a domain graph is generated from a system of trajectories.
For the applications discussed in this paper we are interested in particular trajectories that can be defined in terms of the domain graph.

\begin{defn}
\label{defn:domainwall}
Let $D = (V, E)$ be the domain graph generated by $\mathcal{ST}$ on $\cX$.
A trajectory which is the finite concatenation of wall trajectories is said to be a \emph{domain-wall trajectory.}
The \emph{associated domain graph path} of a domain-wall trajectory $x$ is the path of the domain graph edges corresponding to the wall trajectories which comprise $x$.
\end{defn}

\begin{defn}
\label{defn:searchgraph}
Let $\mathcal{ST}$ be a system of trajectories on $X = (0,\infty)^N$ with rectangular decomposition $\cX$.
If every domain trajectory is monotonic and every wall trajectory undergoes at most one extremal event, we say $\mathcal{ST}$ is \emph{extrema-pattern-matchable} with respect to $\cX$.
In this case, the labeled directed graph $\cS = (V, E, \Sigma_{\text{ext}}^N, \ell)$ is said to be a \emph{search graph} if $(V,E)$ is the domain graph and $\ell$ reflects, as follows, our level of knowledge of the behaviors of trajectories:

\[ \ell(\xi')_n = \begin{cases}

I \mbox{ if we know $x_n(t)$ is increasing for every trajectory $x(t)$ in domain $\xi'$, else}\\
D \mbox{ if we know $x_n(t)$ is decreasing for every trajectory $x(t)$ in domain $\xi'$, else}\\
* \mbox{ otherwise }
 \end{cases}
 \]

 \[ \ell(\xi \to \xi')_n = \begin{cases}
 - \mbox{ if we've ruled out local extrema for $x_n$ on the wall between $\xi$ and $\xi'$, else}\\
m \mbox{ if we've ruled out local maxima for $x_n$ on the wall between $\xi$ and $\xi'$, else}\\
M \mbox{ if we've ruled out local minima for $x_n$ on the wall between $\xi$ and $\xi'$, else}\\
* \mbox{ otherwise}
 \end{cases}
 \]

\end{defn}

Observe that $*$ indicates a lack of knowledge.
If a  system of trajectories is extrema-pattern-matchable, we can always make its domain graph into a search graph by choosing all labels to be $*$.
This would lead to a higher rate of false positives in matching; it is better to assign the strongest labels one can prove.

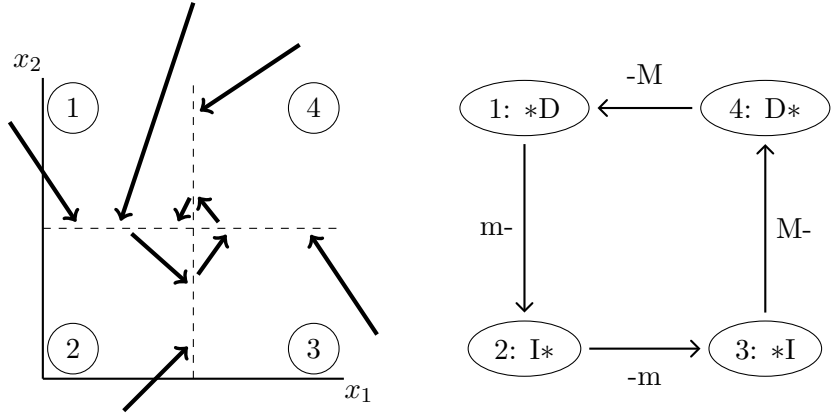
\begin{figure}
\centering
\begin{tabular}{cc}
	\begin{tikzpicture}
		\draw[thick]{(0,0) -- (0,-4)};
		\draw[thick]{(0,-4) -- (4,-4)};
		\draw[dashed]{(2,-4) -- (2,0)};
		\draw[dashed]{(0,-2) -- (4,-2)};
		\draw[->,shorten >=3pt,ultra thick]{(2,1) -> (1,-2)};
		\draw[->,shorten >=3pt,shorten <=3pt,ultra thick]{(1.1,-2) -> (2,-2.8)};
		\draw[->,shorten >=3pt,shorten <=3pt,ultra thick]{(2,-2.7) -> (2.5,-2)};
		\draw[->,shorten >=3pt,shorten <=3pt,ultra thick]{(2.4,-2) -> (2,-1.5)};
		\draw[->,shorten >=3pt,shorten <=3pt,ultra thick]{(2,-1.5) -> (1.75,-2)};

		\draw[->,shorten >=3pt,shorten <=3pt,ultra thick]{(-0.5,-0.5) -> (0.5,-2)};
		\draw[->,shorten >=3pt,shorten <=3pt,ultra thick]{(1,-4.5) -> (2,-3.5)};
		\draw[->,shorten >=3pt,shorten <=3pt,ultra thick]{(4.5,-3.5) -> (3.5,-2)};
		\draw[->,shorten >=3pt,shorten <=3pt,ultra thick]{(3.5,0.5) -> (2,-0.5)};

      \node () at (4.2,-4.2) {$x_1$};
      \node () at (-0.2,0.2) {$x_2$};
      
      \node[style={circle,draw}] () at (0.4,-0.4) {1};
      \node[style={circle,draw}] () at (0.4,-3.6) {2};
      \node[style={circle,draw}] () at (3.6,-3.6) {3};
      \node[style={circle,draw}] () at (3.6,-0.4) {4};

      \node[style={ellipse,draw}] (1) at (6.4,-0.4) {1: $*$D};
      \node[style={ellipse,draw}] (2) at (6.4,-3.6) {2: I$*$};
      \node[style={ellipse,draw}] (3) at (9.6,-3.6) {3: $*$I};
      \node[style={ellipse,draw}] (4) at (9.6,-0.4) {4: D$*$};

      \path[->,>=angle 90,thick]
      (1) edge[shorten <= 3pt, shorten >= 3pt,left] node[] {m-} (2)
      (2) edge[shorten <= 3pt, shorten >= 3pt] node[below=5pt] {-m} (3)
      (3) edge[shorten <= 3pt, shorten >= 3pt,right] node[] {M-} (4)
      (4) edge[shorten <= 3pt, shorten >= 3pt] node[above=5pt] {-M} (1)
      ;

	\end{tikzpicture}
\end{tabular}
\caption{Left: An example system of trajectories that is extrema-pattern-matchable over a rectangular decomposition of four components. Right: The associated search graph.}\label{fig:exampleST}
\end{figure}

As an example, consider a system of trajectories over a rectangular decomposition of $\R^2$, with trajectories  qualitatively depicted in Figure~\ref{fig:exampleST} (left).
The walls are shown as  dotted lines. The domains are labeled 1-4, and within each domain the trajectories are monotonic in each variable, satisfying the requirements to be  extrema-pattern-matchable. Within domain 1, $x_2$ trajectories are monotonically decreasing, while $x_1$ trajectories either monotonically decrease or monotonically increase. The associated search graph shown in Figure~\ref{fig:exampleST} (right); node 1 corresponding to domain 1 is labeled $*$D and  likewise for the other nodes. The edges  between nodes in the search graph correspond to concatenation of domain trajectories. Clearly, there cannot be a maximum in $x_1$ as we pass from domain 1 to domain 2, but there could be a minimum, whereas $x_2$ is constantly decreasing and cannot have either a minimum or maximum on that wall. The edge (1 $\to$ 2) in the search graph is therefore labeled m-, and similar arguments hold for the other edges. 

%!TEX root = ../PatternMatchPaper.tex

\subsection{Matching Pattern Graphs against Search Graphs}
\label{section:matchingtheorem}

As indicated in the introduction we are interested in identifying whether our model for dynamics is capable of producing sequences of maxima and minima that do not contradict the experimental data.
The capability is equivalent to the existence of a matching between a pattern graph and the search graph.
To do this we impose a particular matching relation.
\begin{defn}
\label{defn:extremalmatching}
The \emph{extremal event matching relation}  $\sim_{\text{ext}}$ on $\Sigma_{\text{ext}}$ is given by
\begin{enumerate}
\item (Vertices) $(I \sim_{\text{ext}} *)$, $(I \sim_{\text{ext}} I)$, $(D \sim_{\text{ext}} *)$, $(D \sim_{\text{ext}} D)$, $(* \sim_{\text{ext}} *)$,
\item (Edges) $(- \sim_{\text{ext}} -)$, $(- \sim_{\text{ext}} m)$, $(- \sim_{\text{ext}} M)$, $(- \sim_{\text{ext}} *)$, $(m \sim_{\text{ext}} m)$, $(m \sim_{\text{ext}} *)$, $(M \sim_{\text{ext}} M)$, $(M \sim_{\text{ext}} *)$, $(* \sim_{\text{ext}} *)$.
\end{enumerate}

Given a pattern graph $\cP$ and a search graph $\cS$, we extend this relation to $\Sigma_{\text{ext}}^*$ by defining $a \sim_{\text{ext}} b$ whenever for all $1 \leq i \leq N$, $a_i \sim_{\text{ext}} b_i$.
\end{defn}

\begin{thm}
\label{thm:NoFalseNegatives}
Let $\cP$ be a pattern graph for a poset of extrema $(P,<,\mu)$ and let $\cS$ be a search graph for a  system of trajectories $\mathcal{ST}$ which is extrema-pattern-matchable with respect to a cubical decomposition $\cX$ of $X = (0,\infty)^N$. 
If $\mathcal{ST}$ admits a domain-wall trajectory with a sequence of extremal events corresponding to a linear extension of $P$, then there exists a path $p \in \cP$ from root to leaf, and a path $s$ in $\cS$ such that $p \sim_{\text{ext}} s$.
\end{thm}

%!TEX root = ../PatternMatchPaper.tex

%\subsection{Proof of Theorem~\ref{thm:NoFalseNegatives}}
%\label{section:NoFalseNegativesProof}

\begin{proof}
Let $<'$ be a linear extension of $P$ and name the elements of $P$ as $e_1 <' e_{2} <' \cdots <' e_n$.
Suppose that $\phi : [t_0, t_1] \to X$ is a domain-wall trajectory with the sequence of extremal events $e_1, e_2, \cdots, e_n$.
We show there exists a path $p$ from root to leaf in $\cP$ and a path $s$ in $\cS$ such that $p \sim s$.

\emph{Step 1.}
We construct a path $p$ from root to leaf in $\cP$ and a path $s$ in $\cS$.
By Definition~\ref{defn:domainwall}, since $\phi$ is a domain-wall trajectory, it can be written as a concatenation of wall trajectories $\phi^i : [t_i, t_{i+1}]$ for $i = 1, 2, \cdots, m$.
Because $\cS$ is a search graph, Definition~\ref{defn:searchgraph} implies that extremal events for $\phi(t)$ can only occur on walls (i.e. during times when $\phi(t) \in X_{d-1}$) and at most one kind of extremal event can occur on a given wall.
Since it is impossible for the same extremal event to occur twice in a row (e.g. between any two local minima there must be an intervening local maxima), and wall trajectories intersect precisely one wall, it follows that each wall trajectory experiences at most one extremal event. If we denote the set of extremal events which occur on the wall trajectory $\phi^i$ as $E_i$, then $\mbox{card } E_i \leq 1$. Therefore, $E_i = \emptyset$ or $\{e_j\}$ for some $j$.
Since $\phi$ experiences the events $e_1, e_2, \cdots, e_n$ in order, and $\phi$ is the concatentation of the trajectories $\phi^i$, it follows that there exists an increasing function $\mu : \{1,\cdots,n\} \to \{1, \cdots, m\}$ such that for $i \in \{1, \cdots, n\}$, $e_i \in E_{\mu(i)}$.
Define $p_1 := \emptyset$, and for $i \in \{1, \cdots, m\}$ define 
$p_{i+1} := \bigcup_{j=1}^{i} E_j$.

We show $p_1 \to p_2 \to \cdots \to p_{m+1}$ is a path in $\cP$ from root to leaf.
To this end it suffices to show that: (1) for each $i \in \{1,\cdots,m+1\}$, $p_i$ is a down set of $P$, (2) for each $i \in \{1, \cdots, m\}$ there is an edge $p_i \to p_{i+1}$ in $\cP$, 
(3) $p_1 = \emptyset$, and (4) $p_{m+1} = P$.

Let $i \in \{1,\cdots,m+1\}$.
Define $k = \max \mu^{-1}(\{1,\cdots,i\})$.
Since $\mu$ is increasing it follows that 
$p_{i+1} = \bigcup_{j=1}^{i} E_j = \{e_1, e_2, \cdots, e_k\}$.

Since $e_1 <' e_{2} <' \cdots <' e_{n}$ and $<'$ is a linear extension of $P$, it follows that $p_i$ is a down set of $P$. 
This demonstrates (1).
Now let $i \in \{1, \cdots, m\}$.
We show $p_i \to p_{i+1}$ is an edge in $\cP$.
There are two cases: either (a) $E_i = \emptyset$ and $p_i = p_{i+1}$, or else (b) or else $E_i = \{e_k\}$ for some $k$ and $p_{i+1} = p_{i} \cup \{e_k\}$.
For case (a), $p_i \to p_{i+1}$ is an edge in $\cP$ since the pattern graph admits all self edges.
For case (b), $p_i \to p_{i+1}$ an edge in $\cP$ since $\cP$ contains the edges present in the down set graph of $P$.
This demonstrates (2).
That $p_1 = \emptyset$ is by definition. 
This demonstrates (3).
Finally,  $p_{m+1} = \bigcup_{i=1}^m E_i = P$.
This demonstrates (4).
Since (1), (2), (3), and (4) hold we have that $p = p_1 \to p_2 \to \cdots \to p_{m+1}$ is a path from root to leaf in $\cP$.
Let $s$ be the path in $\cS$ corresponding to the sequence of wall trajectories $\phi^i$ (i.e. the path associated with the domain-wall trajectory $\phi$).
Denote the vertices of the path $s$ in order as $s_1 \to s_2 \to \cdots \to s_{m+1}$.
Note that the wall trajectories $\phi^i$ correspond to the edges $s_i \to s_{i+1}$ in $s$.
We have constructed a path $p$ from root to leaf in $\cP$ and a path $s$ in $\cS$, completing Step 1.

\emph{Step 2.}
We show that for $p$ and $s$ so constructed, $p \sim s$ holds.
By Definition~\ref{defn:MatchingPaths} and Definition~\ref{defn:extremalmatching}, it suffices to show that for each $i \in \{1,\cdots, N\}$, for each $j \in \{1, \cdots, m+1 \}$, $\ell(p_j)_i \sim' \ell(s_j)_i$   (i.e. vertex labels match) and for each  $i \in \{1,\cdots, N\}$  for each $j \in \{1, \cdots, m \}$, $\ell(p_j \to p_{j+1})_i \sim' \ell(s_j \to s_{j+1})_i$ (i.e. edge labels match).

\emph{Proof that edge labels match:}
Let $i \in \{1,\cdots, N\}$  and $j \in \{1,\cdots,m\}$.
We show
\begin{equation}
\label{eqn:edgelabelmatch}
\ell(p_j \to p_{j+1})_i \sim' \ell(s_j \to s_{j+1})_i.
\end{equation}
There are two cases: either (1) $E_j = \emptyset$, or (2) $E_j = \{e_k\}$ for some $k$.
For case (1), $E_j = \emptyset$ implies $p_j = p_{j+1}$ and hence $\ell(p_j \to p_{j+1})_i = -$.
Meanwhile $\ell(s_j \to s_{j+1})_i \in \{ -, m, M, *\}$.
By Definition~\ref{defn:extremalmatching} it follows that Equation~\eqref{eqn:edgelabelmatch} holds for case (1).
For case (2), $E_j = \{e_k\}$ for some $k$, we distinguish three subcases: (a) $e_k$ is local minimum for variable $i$, (b) $e_k$ is a local maximum for variable $i$, or (c) $e_k$ is not a local extremum for variable $i$.
For subcase (a), $\ell(p_j \to p_{j+1})_i = m$.
Since the wall trajectory $\phi^j$ experienced a local minimum for variable $i$, it follows that we could not have ruled out a local minimum on the wall corresponding to the edge $s_j \to s_{j+1}$.
This eliminates the possibility that $\ell(s_j \to s_{j+1})_i$ is either $-$ or $M$, i.e. $\ell(s_j \to s_{j+1})_i \in \{m, *\}$. Since $m \sim' m$ and $m \sim' *$, Equation~\eqref{eqn:edgelabelmatch} holds for subcase (a).
Subcase (b) is similar.
For subcase (c), $\ell(p_j \to p_{j+1})_i = -$, and the argument of case (1) again applies.
Hence Equation~\eqref{eqn:edgelabelmatch} holds in all cases.

\emph{Proof that vertex labels match:}
Let $i \in \{1,\cdots, N\}$  and $j \in \{1,\cdots,m+1\}$.
We show $\ell(p_j)_i \sim' \ell(s_j)_i$.

Let $P_i = \{ p \in P : \ell(p)_i \in \{I,D\} \}$.
We consider two cases: (1) $P_i = \emptyset$, and (2) $P_i \neq \emptyset$. For case (1), by Definition~\ref{defn:patterngraph}, $P_i = \emptyset$ implies $\ell(p_j)_i = *$.
By Definition~\ref{defn:searchgraph}, $\ell(s_j)_i \in \{I, D, *\}$, and by Definition~\ref{defn:extremalmatching} $* \sim' I$, $* \sim' D$, and $* \sim' *$.
It follows that $\ell(p_j)_i \sim' \ell(s_j)_i$ for case (1).
For case (2), we assume $P_i \neq \emptyset$.
Then $\ell(p_j)_i \in \{I, D\}$.
There are four subcases depending on whether (a) $\ell(p_j)_i = I$ or $\ell(p_j)_i = D$, and (b) whether $p_j \cap P_i = \emptyset$.
As they are all similar, we only consider the subcase when $\ell(p_j)_i = I$ and and $p_j \cap P_i \neq \emptyset$.
Let $\phi' : [t_1, t_{j}] \to X$ be the domain-wall trajectory $\phi'$ obtained by concatenating $\phi^1$, $\phi^{2}$, $\cdots$, $\phi^{j-1}$.
By construction, $p_j$ is the set of events in $P$ which occur on $\phi'$.
Let $e$ be the maximal element of $p_j \cap P_i$.
By Definition~\ref{defn:patterngraph}, $\ell(p_j)_i = I$ implies that $e$ is a local minimum.
It follows that $\phi'$ is increasing in variable $i$ after event $e$ occurs.
This implies that for sufficiently small $\epsilon > 0$, $\phi^{j-1}\vert_{[t_{j} - \epsilon, t_j]}$ is an increasing trajectory with image contained in the domain $s_j$.
By Definition~\ref{defn:se{}archgraph}, it follows that $\ell(s_j)_i \in \{I, *\}$.
By Definition~\ref{defn:extremalmatching}, $I \sim' I$ and $I \sim' *$, and $\ell(p_j)_i \sim' \ell(s_j)_i$ follows.
Similar arguments for the other three subcases show $\ell(p_j)_i \sim' \ell(s_j)_i$ for case (2).
We have shown $p \sim s$, which completes Step 2.

Since in Step 1 we constructed a path $p$ in $\cP$ from root to leaf and a path $s$ in $\cS$ and in Step 2 we showed $p \sim s$, the proof is complete.
\end{proof}

To continue our example, we take the pattern graph $\cP$ from Figure~\ref{fig:earlypoexample} (right) and the search graph $\cS$ from Figure~\ref{fig:exampleST} (right) and seek matching paths $p \in \cP$, $s \in \cS$. To do this, we form the alignment graph as in Definition~\ref{defn:alignmentgraph} using the matching relation $\sim_{\text{ext}}$ given in Definition~\ref{defn:extremalmatching}. We then apply Proposition~\ref{prop:AlignmentGraphProperty} that states that finding paths in the alignment graph is equivalent to finding pairs of matching paths in the pattern and search graphs. In particular, we seek a match to a path $p \in \cP$ that is a linear extension of the poset of extrema in Figure~\ref{fig:earlypoexample} (left), to verify that the system of trajectories $\mathcal{ST}$ in Figure~\ref{fig:exampleST} (left) can support the constraints on the order of extrema summarized by the poset.

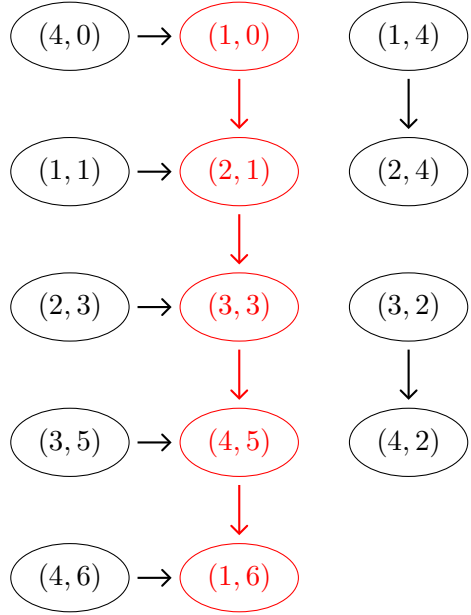
\begin{figure}
   \begin{center}
    \begin{tikzpicture}[main node/.style={ellipse,fill=white!20,draw},scale=1.8]
      \node[main node] (40) at (-1,-1) {$(4,0)$};
      \node[main node,draw=red,text=red] (10) at (0.25,-1) {$(1,0)$};
      \node[main node] (11) at (-1,-2) {$(1,1)$};
      \node[main node] (14) at (1.5,-1) {$(1,4)$};

      \node[main node,draw=red,text=red] (21) at (0.25,-2) {$(2,1)$};
      \node[main node] (23) at (-1,-3) {$(2,3)$};
      \node[main node] (24) at (1.5,-2) {$(2,4)$};

      \node[main node,draw=red,text=red] (33) at (0.25,-3) {$(3,3)$};
      \node[main node] (35) at (-1,-4) {$(3,5)$};
      \node[main node] (32) at (1.5,-3) {$(3,2)$};

      \node[main node] (46) at (-1,-5) {$(4,6)$};
      \node[main node,draw=red,text=red] (45) at (0.25,-4) {$(4,5)$};
      \node[main node] (42) at (1.5,-4) {$(4,2)$};

      \node[main node,draw=red,text=red] (16) at (0.25,-5) {$(1,6)$};

      \path[->,>=angle 90,thick]
      (40) edge[shorten <= 3pt, shorten >= 3pt] node[] {} (10)
      (10) edge[shorten <= 3pt, shorten >= 3pt,color=red] node[] {} (21)
      (11) edge[shorten <= 3pt, shorten >= 3pt] node[] {} (21)
      (14) edge[shorten <= 3pt, shorten >= 3pt] node[] {} (24)

      (21) edge[shorten <= 3pt, shorten >= 3pt,color=red] node[] {} (33)
      (23) edge[shorten <= 3pt, shorten >= 3pt] node[] {} (33)

      (33) edge[shorten <= 3pt, shorten >= 3pt,color=red] node[] {} (45)
      (35) edge[shorten <= 3pt, shorten >= 3pt] node[] {} (45)
      (32) edge[shorten <= 3pt, shorten >= 3pt] node[] {} (42)

      (45) edge[shorten <= 3pt, shorten >= 3pt,color=red] node[] {} (16)
      (46) edge[shorten <= 3pt, shorten >= 3pt] node[] {} (16)

      ;

      \end{tikzpicture}
    \end{center}\caption{Alignment graph for the pattern graph in Figure~\ref{fig:earlypoexample} (right) and the search graph in Figure~\ref{fig:exampleST} (right). The red path indicates a match between paths in the graphs.}\label{fig:examplealigngraph}
\end{figure}

The alignment graph is given in Figure~\ref{fig:examplealigngraph}, where each node is labeled by a pair $(a,b)$ where $a$ is a node identifier for the search graph (integers 1-4) and $b$ is a node identifier for the pattern graph (integers 0-6). The red path denotes a match between path  $p = (0,1,3,5,6)$ in the pattern graph  in Figure~\ref{fig:earlypoexample} (right) and cyclic path $s = (1,2,3,4,1)$ in the search graph in Figure~\ref{fig:exampleST} (right). We notice that $p = (0,1,3,5,6)$ corresponds to a linear extension of the poset of extrema in Figure~\ref{fig:earlypoexample} (left), since it is a path from root to leaf of the pattern graph (Theorem~\ref{thm:bijection}). Therefore $\mathcal{ST}$ has at least one trajectory with a sequence of extrema respecting the constraints of the poset of extrema.

\section{Application to Regulatory Networks}
\label{section:Application}
%!TEX root = ../PatternMatchPaper.tex

%\subsection{Overview}
%\label{section:ApplicationOverview}

As indicated in the introduction to provide a demonstration of how to apply the combinatorial tools described in the previous sections we make use of DSGRN.  A complete description of the mathematical framework can be found in \cite{cummins2016combinatorial}, however for the benefit of the reader we begin this section with a short review.
We then present an application to simple system associated with the cell cycle of \textit{S. cerevisiae} using experimental time-series data (provided courtesy of the Haase lab; see~\cite{leman2014analyzing} for data collection methods) for mRNA sequences associated with SWI4, HCM1, NDD1, and YOX1 collected at time intervals of 5 minutes.

%!TEX root = ../PatternMatchPaper.tex

\subsection{DSGRN Model for Regulatory Networks}
\label{section:DSGRN}
We provide a mathematical definition of a regulatory network and its associated parameters.
We use this to construct a system of trajectories and prove some simple results concerning its structure.
We conclude by relating the system of trajectories to the output of DSGRN which provides us with a means of  analyzing specific data sets.
 
\begin{defn}
\label{defn:rn}
A \textit{regulatory network} $\mathbf{RN}=(V, E, M)$ consists of vertices $V = \setof{1,\ldots,N}$ called \textit{network nodes}, annotated directed edges $E \subset V\times V \times \{ \to, \dashv \} $ called \textit{interactions}, and for each $k \in V$, polynomial monotone increasing functions $M_n : \R^{|S_n|} \to \R$ called \emph{node logics} where  $S_n := \{ (i,n,\cdot) \in E \}$ is called the \emph{source set} of $k$. 

An $\to$ annotated edge is referred to as an \textit{activation} and an $\dashv$ annotated edge is called a \textit{repression}.
We indicate that either $i \to j$ or $i \dashv j$ without specifying which by writing $(i, j) \in E$.
We allow self-edges. 
From one node to another we admit at most one type of annotated edge, e.g.\ we cannot have both $i \to j$ and $i \dashv j$ simultaneously.
The $n$-th \emph{target set} is given by $T_n := \{ (n,j) \in E  \}$.
\end{defn}

A parameterized family of dynamics is generated from the regulatory network. 

\begin{defn}
\label{defn:paramspace}
A \emph{parameter} for a regulatory network $\mathbf{RN} = (V,E,M)$ is a tuple $z \in Z \subset (0,\infty)^{(N+3\cdot |E|)}$.
The coordinates of a parameter $z$ are associated with the nodes and edges of $\mathbf{RN}$ and are given by the values of four functions $\gamma : V \to (0,\infty)$, and $\ell,u,\Theta : E \to (0, \infty)$ with the constraint that $\ell(e) \leq u(e)$ for each $e\in E$.
\end{defn}

The functions $\gamma$, $\ell$, $u$, and $\Theta$ are used to decompose phase space and generate dynamics as follows.
Define 
\[
\Theta_n := \setof{\Theta((n,j)) : (n,j) \in T_n},\quad \text{for $n \in \{1,\cdots,N\}$}.
\]
and assume that for all $n = 1,\ldots, N$,
\begin{equation}
\label{eq:thetanotequal}
\text{if $\Theta((n,j)), \Theta((n,k))\in \Theta_n$, then $\Theta((n,j))  \neq \Theta((n,k))$.}
\end{equation}
Then, $(\Theta_1, \Theta_2, \cdots, \Theta_N)$ defines a rectangular decomposition $\cX$ (see Definition~\ref{defn:cubicaldecomposition}) on $X := (0,\infty)^N$.

Define $\Gamma$ to be the diagonal $N \times N$ matrix with diagonal entries $\gamma(n)$ for $n \in \{1,\cdots,N\}$.
Define $ W : E \times X \to (0, \infty)$ via 
\[ 
W((i,j),x) = \begin{cases} \ell((i,j)) & \text{if $x_i < \Theta((i,j))$ and $x_i \to x_j $,}\\
\ell((i,j)) & \text{if $ x_i > \Theta((i,j))$ and $x_i \dashv x_j$,} \\
u((i,j)) & \text{if $ x_i > \Theta((i,j))$ and $ x_i \to x_j$,} \\
u((i,j)) & \text{if $ x_i < \Theta((i,j))$ and $ x_i \dashv x_j$,} \\
0 & \text{otherwise.}
 \end{cases}.
 \]
Finally, define  $\Lambda : X \to \R^N$ by 
\[ 
\Lambda_n(x) := M_n \circ W\vert_{S_n \times X}.
\]
We are interested in dynamics generated by differential equations of the form
\[
\dot{x} = -\Gamma x + \Lambda(x).
\]
Observe that if $\xi\in\cX_N$, then $\Lambda$ is constant on $\xi$ and hence it makes sense to write $\Lambda(\xi)$.

\begin{defn}
A parameter value $z\in Z$ is \emph{regular} if \eqref{eq:thetanotequal} holds,  $\ell(e) < u(e) \quad \forall e \in E$,  and 
\begin{equation}
-\gamma(n)\Theta(n,k) + \Lambda_n(\xi) \neq 0
\end{equation}
if an $N-1$ dimensional face of  $\xi\in\cX_N$ lies in the hyperplane defined by $x_n = \Theta(n,k)$.
The set of regular parameter values is denoted by $Z^R$.
\end{defn}

\begin{defn}
\label{defn:DSGRNDynamicalSystem}
A $\mathbf{RN}$ \emph{domain trajectory} is a function $x\colon [t_0,t_1] \to \cl(\xi)$, where $\xi\in\cX_N$, that solves the  differential equation 
\begin{equation}
\label{eq:switchingsystemODE}
\dot{x} = -\Gamma x + \Lambda(\xi).
\end{equation}
Given  $z\in Z^R$, the associated $\mathbf{RN}$ \emph{system of trajectories} at parameter value $z$, denoted $\mathcal{ST}(\mathbf{RN},z)$, is defined to be the smallest system of trajectories (see Definition~\ref{defn:dynamicalsystem}) which contains every $\mathbf{RN}$ domain trajectory.
\end{defn}

\begin{rem}
It is straightforward to verify that under Definition~\ref{defn:dynamicalsystem}, the intersection of two systems of trajectories is again a system of trajectories.
Thus the notion of the smallest  system of trajectories containing some set of trajectories is well-defined.
\end{rem}

Observe that \eqref{eq:switchingsystemODE} is a linear differential equation and that
\[
P^\xi := \Gamma^{-1}\Lambda(\xi) \in \R^N
\]
is a globally attracting fixed point within $\cl(\xi)$.

Let $\pi_n\colon \R^N \to \R$ be the canonical projection map onto the $n$-th coordinate.

\begin{prop}
\label{prop:dsgrnedges}
Let $\mathbf{RN}$ be a regulatory network and $z \in Z^R$.
Consider the system of trajectories $\mathcal{ST}(\mathbf{RN},z)$.
Let $\xi, \xi' \in \cX_N$ be  separated by the hyperplane $x_n = \Theta((n,j))$ for some $(n,j) \in E$ such that $\pi_n(\xi) < \pi_n(\xi')$.
Then,
 \begin{enumerate}
 \item[(i)] there exists a wall trajectory from $\xi$ to $\xi'$ if and only if $\max\{P^\xi_n, P^{\xi'}_n\} > \Theta((n,j))$
 \item[(ii)] there exists a wall trajectory from $\xi'$ to $\xi$ if and only if $\min\{P^\xi_n, P^{\xi'}_n\} < \Theta((n,j))$
 \end{enumerate}
 \end{prop}
\begin{proof}
We show (i) and leave (ii) to the reader.
Suppose there exists a wall trajectory $x : [0, T] \to (0,\infty)^N$ from $\xi$ to $\xi'$.
By Definition~\ref{defn:dynamicalsystem}, the restrictions $x\vert_A$ and $x\vert_B$ onto $A = x^{-1}(\bar{\xi})$ and $B = x^{-1}(\bar{\xi'})$ are again trajectories.
By Definition~\ref{defn:DSGRNDynamicalSystem}, $x\vert_A$ and $x\vert_B$  are solutions to Equation~\eqref{eq:switchingsystemODE}.
Such solutions are monotonic, so it follows that $x\vert_A$ and $x\vert_B$ are increasing.
This requires $\pi_n(P^\xi) > \Theta((n,j))$ and $\pi_n(P^{\xi'}) > \Theta((n,j))$ (with strictness since we reach or leave the wall in finite time), yielding $\max\setof{P^\xi_n, P^{\xi'}_n} > \Theta((n,j))$ as desired.

To prove the converse suppose $\max\setof{P^\xi_n, P^{\xi'}_n} > \Theta((n,j))$.
Let $\hat{x} \in \sigma$, where $\sigma\in\cX_{N-1}$ is the cell between $\xi$ and $\xi'$.
Solve the initial value problem \eqref{eq:switchingsystemODE} with initial value $\hat{x}$ in forward time in $\xi'$ and in backward time in $\xi$ to obtain solutions $x \colon [t_0, t_1]\to \cl(\xi)$ and $y \colon [t_1, t_2]\to \cl(\xi')$ such that $x(t_1) = y(t_1) = \hat{x}$.
By Definition~\ref{defn:DSGRNDynamicalSystem}, $x$ and $y$ are trajectories in $\mathcal{ST}(\mathbf{RN},z)$.
By Definition~\ref{defn:dynamicalsystem} the concatenation of $x$ and $y$ is again a trajectory.
This yields a wall trajectory from $\xi$ to $\xi'$.  
\end{proof}

\begin{prop}
\label{prop:DSGRNSearchGraphProp}
Let $\mathbf{RN}$ be a regulatory network,  $z \in Z^R$, and $\xi\in\cX_N$ be a domain.
Then $\mathcal{ST}(\mathbf{RN},z)$ has the following properties:
\begin{enumerate}
\item[(i)] Every trajectory $x(t)$ in  $\xi$ is monotonic in each variable.
\item[(ii)] If $P^\xi_n > \pi_n(\xi)$, then for every trajectory $x(t)$ in $\xi$, $x_n(t)$ is an increasing function.
\item[(iii)] If $P^\xi_n < \pi_n(\xi)$, then for every trajectory $x(t)$ in $\xi$, $x_n(t)$ is a decreasing function.
\item[(iv)] If $P^\xi_n \in  \pi_n(\xi)$, there exist trajectories $x(t)$ in $\xi$ where $x_n(t)$ may be either an increasing, decreasing, or constant function.
\item[(v)] Let $w$ be a wall associated with the hyperplane $x_n = \Theta((n,j))$ arising from the regulatory network interaction $x_n \to x_j$. 
Then, the only type of extremum a wall trajectory can undergo as it passes through $w$ is a local minimum in the variable $x_j$.
\item[(vi)] Let $w$ be a wall associated with the hyperplane $x_n = \Theta((n,j))$ arising from the regulatory network interaction $x_n \dashv x_j$. 
Then, the only type of extremum a wall trajectory can undergo as it passes through $w$ is a local maximum in the variable $x_j$.
\end{enumerate}
\end{prop}
\begin{proof}
(i)-(iv) follow immediately from Equation~\eqref{eq:switchingsystemODE}.

We show (v) and leave (vi) to the reader.
Let $w$ be a wall associated with the hyperplane $x_n = \Theta((n,j))$ arising from the regulatory network interaction $x_n \to x_j$.
Let $\xi, \xi'$ be the adjacent domains that $w$ separates, such that $\pi_n(\xi) < \pi_n(\xi')$.
Let $x:[t_0,t_1] \to X$ be a wall trajectory from $\xi$ to $\xi'$.
We show that $x$ cannot undergo any kind of extremum except possibly a local minimum in the variable $x_j$.

Since $z\in Z^R$, $\Theta((n,j)) \neq \Theta((n,k))$ for $j \neq k$.
This implies that $P^\xi_k = P^{\xi'}_k$ for all $k \neq j$.
Define  $x\vert_A$ and $x\vert_B$, where $A = x^{-1}(\bar{\xi})$ and $B = x^{-1}(\bar{\xi'})$.
Since $x\vert_A$ and $x\vert_B$ each obey Equation~\eqref{eq:switchingsystemODE} on their respective domains, it follows that $x_k$ obeys $\dot{x}_k = -\gamma_k(x_k - P^u_k)$ everywhere.
Thus $x_k(t)$ is monotonic, and hence experiences no extremal event.
Now we show $x$ cannot undergo a local maximum event in variable $x_j$. Since $l((i,j))<u((i,j))$ it follows from the definitions that we must have $P^\xi_j < P^{\xi'}_j$.
If $x\vert_A$ is constant or decreasing in the $j$th coordinate, then there cannot be a local maximum as we pass the wall.
So we consider only the case where $x\vert_A$ is increasing in the $j$th coordinate.
This case requires that $x_j(A) < P^u_j$.
Hence, $x_j(\min B) < P^\xi_j < P^{\xi'}_j$.
Since $x\vert_B$ is a solution of the initial value problem corresponding to Equation~\eqref{eq:switchingsystemODE} with an initial condition for $x_j$ less than $P^{\xi'}_j$, it follows that $x_j$ is everywhere increasing. 
Therefore, $x_j$ does not experience a local maximum. 
\end{proof}

\begin{thm}
\label{thm:DSGRNYieldsSearchGraphs}
Let $\mathbf{RN}$ be a regulatory network and  $z \in Z^R$.
Let $\mathcal{ST}(\mathbf{RN},z)$ be the associated system of trajectories.
Let $\cS = (V,E,\Sigma,\ell)$ be the labeled directed graph given by
\begin{enumerate}[label=(\roman*)]
\item 
 $V = \cX_N$  
\item $E = \{ (\xi,\xi') \in \cX_N^2 : \xi \text{ and $\xi'$ are adjacent, and for all $n\in V$, either }$ \\
\phantom{indentindentindent} $\left((\pi_n(\xi) < \pi_n(\xi')) \wedge (\min\setof{P^\xi_n,P^{\xi'}_n} > \pi_n(\xi))\right)$ or \\ \phantom{indentindentindent} $\left((\pi_n(\xi) > \pi_n(\xi')) \wedge (\max\setof{P^\xi_n,P^{\xi'}_n} < \pi_n(\xi))\right) \}.$ % \item[(iv)] For all $n\in V$, $\ell(\xi)_n = I$ whenever $P^\xi_n > \pi_n(\xi)$
\item For all $n\in V$, $\ell(\xi)_n = D$ whenever $P^\xi_n < \pi_n(\xi)$
\item For all $n\in V$, $\ell(\xi)_n = *$ whenever $P^\xi_n \in \pi_n(\xi)$
\item For all $n,j,k\in V$, \\ $\ell(\xi \to \xi')_n = -$ whenever $x_j \to x_k$, $n \neq k$, and $x_j = \Theta((j,k))$ separates $\xi $ and $ \xi'$
\item For all $n,j\in V$, \\ $\ell(\xi \to \xi')_n = m$ whenever  $x_j \to x_n$  and $x_j = \Theta((j,n))$ separates  $\xi $ and $ \xi'$ 
\item For all $n,j\in V$, \\$\ell(\xi \to \xi')_n = M$ whenever $x_j \dashv x_n$  and $x_j = \Theta((j,n))$ separates$\xi $ and $ \xi'$.
\end{enumerate}
Then, $\cS$ is a search graph for $\mathcal{ST}(\mathbf{RN},z)$ with the rectangular decomposition $\cX$.
 \end{thm}

\begin{proof}
By Proposition~\ref{prop:dsgrnedges}, it follows that $(V,E)$ is the domain graph for $DSGRN(\mathbf{RN},z)$ with the cubical decomposition $\cX$.
By Proposition~\ref{prop:DSGRNSearchGraphProp}, it follows the vertex and edge labels satisfy the requirements of Definition~\ref{defn:searchgraph}.
\end{proof}

Theorem~\ref{thm:DSGRNYieldsSearchGraphs} guarantees that given a regulatory network and  regular parameter value there exists a search graph.
The next proposition indicates that parameter space admits a finite decomposition, where within each open component of the decomposition the parameters exhibit isomorphic search graphs. 

 \begin{prop}\label{prop:parametergraph}
For a fixed regulatory network the following hold:
\begin{enumerate}
\item[(i)] The regular parameter values $Z^R$ form an open and dense subset of all parameter values $Z$.
\item[(ii)]  $Z^R$ has finitely many connected components.
\item[(iii)]  The connected components of $Z^R$ are semialgebraic sets which can be written as systems of strict inequalities involving polynomials of the parameters.
\item[(iv)]  If $z_1, z_2 \in Z^R$ are in the same connected component of $Z^R$, then the search graph for $\mathcal{ST}(\mathbf{RN},z_1)$ is isomorphic to the search graph associated with $\mathcal{ST}(\mathbf{RN},z_2)$.
\end{enumerate}
\end{prop}

We do not provide a proof of Proposition~\ref{prop:parametergraph} as it is a partial summary of results in \cite{cummins2016combinatorial} that describes the mathematical foundations for the DSGRN software~\cite{DSGRNRepository}.
Given a regulatory network for which $|S_n| \leq 3$ and $|T_n| \leq 3$ the key computational result of \cite{cummins2016combinatorial} is that DSGRN provides an efficient computational scheme for constructing an undirected graph $\mathsf{PG}$, called the parameter graph,  where each node represents one of the connected components described in Proposition~\ref{prop:parametergraph}(iii) and the edges correspond to a notion of adjacency of the parameter regions.
In addition, for each node in the parameter graph DSGRN can be used to compute the associated domain graph, i.e., identify the  set of vertices and the set of edges of $S$ as described in Theorem~\ref{thm:DSGRNYieldsSearchGraphs}(i) and (ii).

From the domain graph, it is possible to extract summary data, called a \emph{Morse Graph}, that provides information about the global dynamics. 
The association of a Morse Graph to each node in the parameter graph $\mathsf{PG}$ gives rise to the notion of a database of dynamical information; the interested reader is referred to \cite{arai2009database, bush2012combinatorial, cummins2016combinatorial} for further details about Morse graphs and dynamical databases. 
For the purposes of this paper, the notion of the domain graph, and the search graph which arises from it, suffices.

 We remark that the system of trajectories $\mathcal{ST}(\mathbf{RN},z)$ qualitatively depicted  in Figure~\ref{fig:exampleST} (left) arises from the regulatory network $\mathbf{RN}(\{x_1,x_2\},\{x_1 \to x_2, x_2 \dashv x_1 \})$ for any regular parameter $z$ satisfying
 \begin{align*}
 \ell((1,2)) &< \theta((1,2)) < u((1,2)) \\
 \ell((2,1)) &< \theta((2,1)) < u((2,1)).
 \end{align*}
%!TEX root = ../PatternMatchPaper.tex

\subsection{Labeled Pattern Graph from Experimental Data}
\label{section:applicationpatterngraph}

\begin{figure}[h!]
\centering
\begin{tabular}{m{3in} m{4in}}
	\hspace{-0.25in}\includegraphics[width=3.5in]{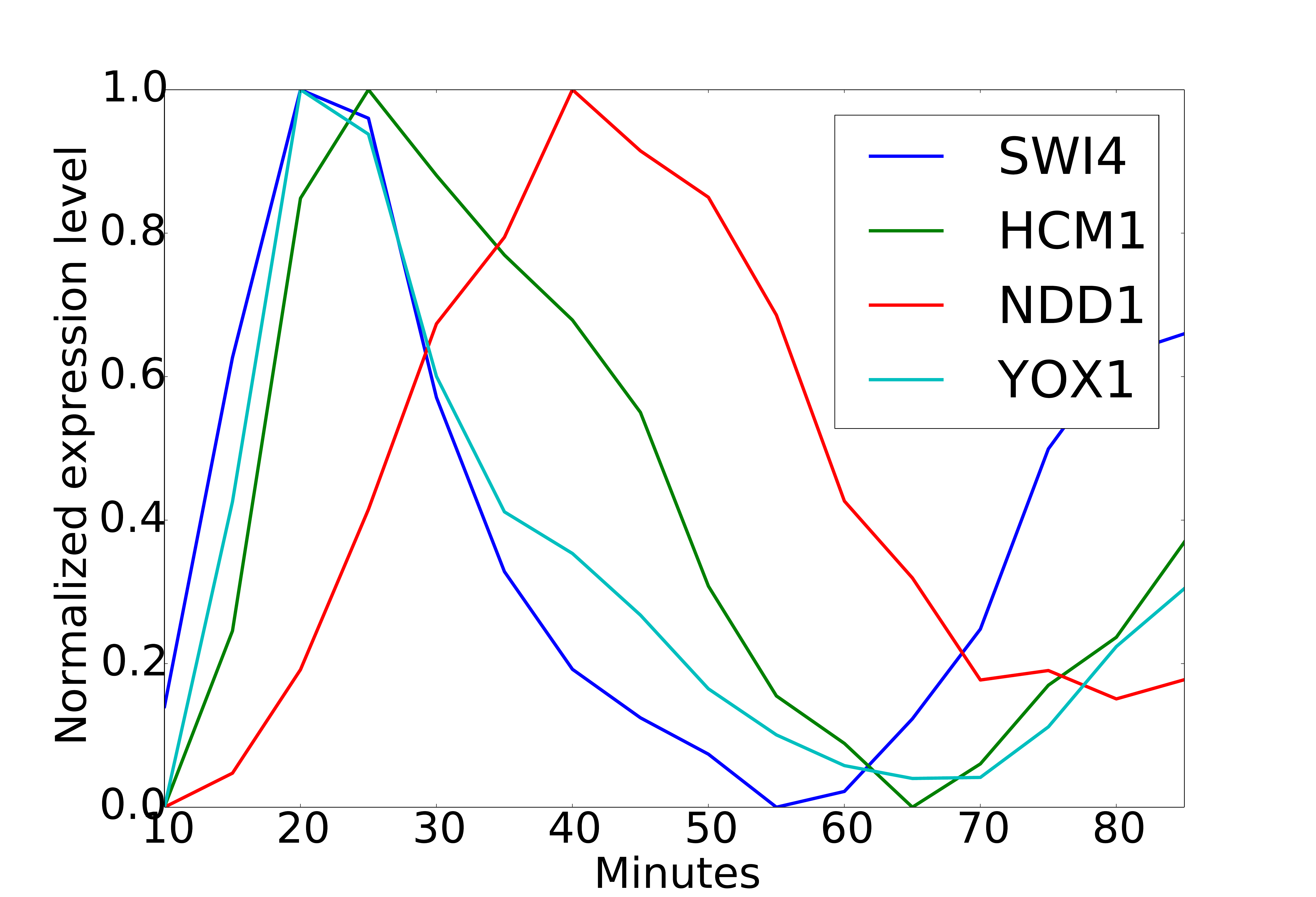}
	&
	 \begin{tabular}{|c|c|c|}
	 \hline
	 event & time interval \\
	 \hline
	 SWI4 min & $(-\infty, 10)$ \\ 
	 HCM1 min & $(-\infty, 10)$ \\ 
	 NDD1 min & $(-\infty, 10)$ \\ 
	 YOX1 min & $(-\infty, 10)$ \\
     
     SWI4 max & $(15,30)$ \\ 
     YOX1 max & $(15,30)$ \\ 
     HCM1 max & $(20,35)$ \\ 
          
     NDD1 max & $(35,45)$\\
          
     SWI4 min & $(50,60)$\\
          
     YOX1 min & $(60,75)$\\ 
     HCM1 min & $(60,75)$\\
          
     NDD1 min & $(70,85)$\\
          
     SWI4 max & $(75,\infty)$\\
          
     HCM1 max & $(85,\infty)$\\
     YOX1 max & $(85,\infty)$\\
          
     NDD1 max & $(85,\infty)$\\

	 \hline
	 \end{tabular}
 \end{tabular}\caption{Left: Time series of  RNA-seq data normalized to range from zero to one. 
  Right: Table of time intervals associated to the extrema in the plot on the left.}
\label{fig:dataset}
\end{figure}

We now turn to the task of generating a labeled pattern graph from experimental data.
The graph in Figure~\ref{fig:dataset}(left) provides normalized expression level data for mRNA sequences associated with SWI4, HCM1, NDD1, and YOX1 from \textit{S. cerevisiae} taken  at time intervals of 5 minutes.
Since we are only concerned with the orderings of the extrema, the normalization of the data makes it easier to identify these extrema.

As indicated in the introduction identifying extrema in data is a serious statistical endeavor that we do not address in this paper.
While our techniques require a set of potential sequences of extrema, they are agnostic with respect to how the potential sequences are derived, therefore we are content for the purpose of this paper to use simple heuristics.
In particular, the table in Figure~\ref{fig:dataset}(right) provides intervals of time within which we declare a maximum or minimum value of expression has occurred.
For example, to allow for noise in the data the tightest time bound we are willing to assume on the maxima for SWI4 and YOX1 is $(15,30)$.
Similarly, we ignore the potential for a local minimum and maximum of NDD1 at time points 70 and 80, and instead assume that a minimum occurs somewhere within the time interval $(70,85)$.

Because we are using intervals to quantify the occurrence in time of extrema we cannot expect to obtain a linear ordering.
Instead we define a partial order $<_\tau$ by 
\begin{equation}
\label{eq:timeIntervalOrdering}
(a,b) <_\tau (c,d) \text{ whenever } b \leq c.
\end{equation}

Note  that the poset of extrema in Figure~\ref{fig:earlypoexample} (left) arises from using $<_\tau$ on rows 1, 4, 5, and 6 in the table; i.e. we form the poset consisting of the first minimum and first maximum of each of $x_1=$ SWI4 and $x_2=$ YOX1.

\begin{figure}[h!]
	\begin{center}
	\includegraphics[width=3.5in]{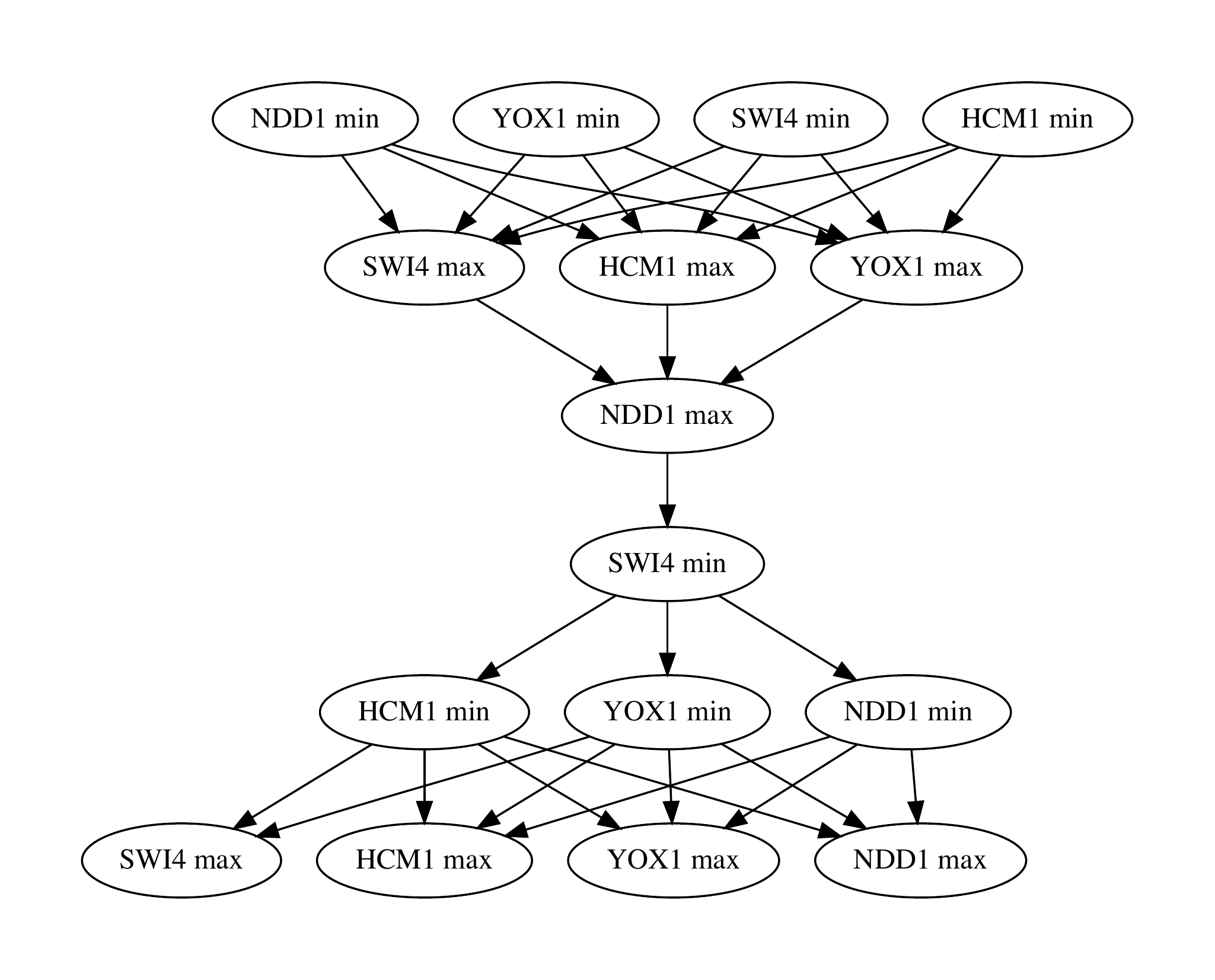}
	\end{center}
	\caption{The pattern (poset) arising from the choice of time intervals of extrema based on the table in Figure~\ref{fig:dataset}. Arrows indicate direction of time.  }
	\label{fig:datapartialorder}
\end{figure}

Using all of the rows in the table in Figure~\ref{fig:dataset} results in the poset indicated in Figure~\ref{fig:datapartialorder}.
Note that the linear extensions of $<_\tau$ correspond to ordered sequences of extrema events.

\begin{figure}
    \begin{center}
	\includegraphics[width=5.5in]{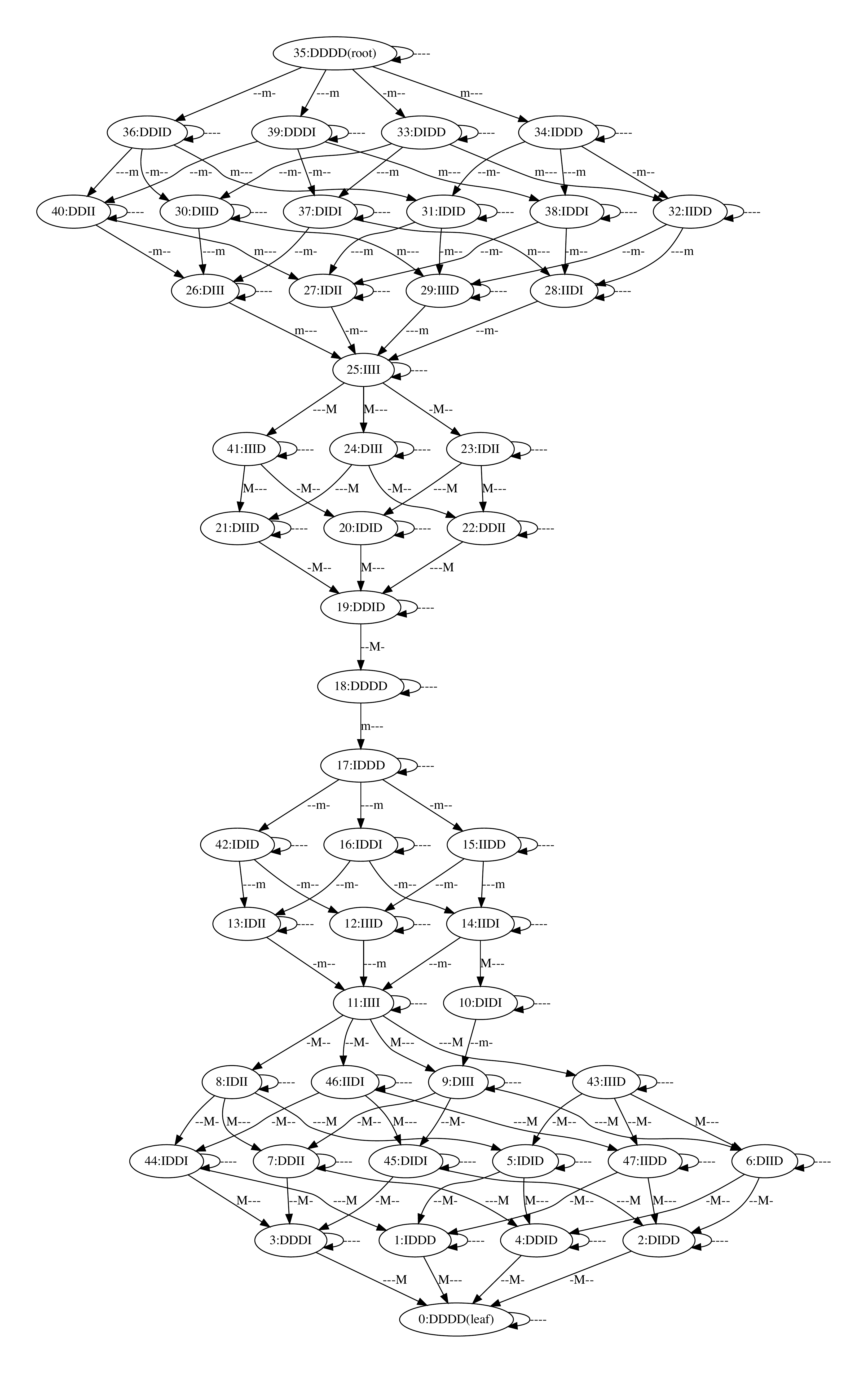}
	\caption{Pattern graph  associated to the pattern in Figure~\ref{fig:datapartialorder}. }
	\label{fig:patterngraph}
	\end{center}
\end{figure}

Observe that we have constructed a poset of extrema $(P,<_\tau;\mu)$ (see Definition~\ref{defn:PosetOfExtrema}) where $P$ consists of the entries of the time interval column in the table in Figure~\ref{fig:dataset} (right), $<_\tau$ is as defined by \eqref{eq:timeIntervalOrdering}, and the values of $\mu \colon P\to \setof{ -,m,M}^4$ are obtained from the event column of the table in Figure~\ref{fig:dataset} (right).
For example, if the first coordinate of $\mu$ corresponds to SWI4, then $P_1 =\setof{ (-\infty,10), (15,30), (50,60), (75,\infty)}$. 
Following Definition~\ref{defn:patterngraph} the associated pattern graph $\cP$ is shown in Figure~\ref{fig:patterngraph}. We remark that Proposition~\ref{prop:downsetcomplexity} applies in this situation, i.e. Algorithm~\ref{alg:posettodownset} can quickly compute the pattern graph $\cP$.

%!TEX root = ../PatternMatchPaper.tex

\subsection{Results for Wavepool Models}
\label{section:ApplicationResults}

\begin{figure}[h!]
  \begin{center}
    \begin{tikzpicture}[main node/.style={rectangle,fill=white!20,draw,font=\sffamily\normalsize},scale=1.6]
      \node[main node] (SWI) at (0,0) {1};
      \node[main node] (HCM) at (1.5,0) {2};
      \node[main node] (NDD) at (3,0) {3};
      \node[main node] (YOX) at (0,1) {4};

      \path[->,>=angle 90,thick]
      (YOX) edge[-|,shorten <= 3pt, shorten >= 3pt, bend left] node[] {} (SWI)
      (SWI) edge[->,shorten <= 3pt, shorten >= 3pt, bend left] node[] {} (YOX)
      (SWI) edge[->,shorten <= 3pt, shorten >= 3pt] node[] {} (HCM)
      (HCM) edge[->,shorten <= 3pt, shorten >= 3pt] node[] {} (NDD)
      (NDD) edge[->,shorten <= 3pt, shorten >= 3pt, bend left] node[] {} (SWI)
      ;

      \end{tikzpicture}
    \end{center}
  \caption{The wavepool regulatory network $\mathbf{RN}_W$  where $M_{\text{1}}$ is multiplication.}\label{fig:regnetwork}
\end{figure}
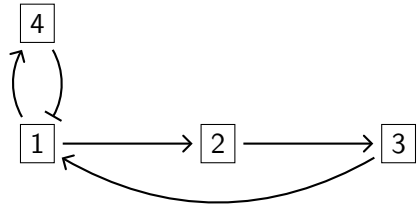

%\begin{figure}[h!]
%  \begin{center}
%    \begin{tikzpicture}[main node/.style={rectangle,fill=white!20,draw,font=\sffamily\normalsize},scale=1.6]
 %     \node[main node] (SWI) at (0,0) {SWI4};
%      \node[main node] (HCM) at (1.5,0) {HCM1};
%      \node[main node] (NDD) at (3,0) {NDD1};
%      \node[main node] (YOX) at (0,1) {YOX1};

%      \path[->,>=angle 90,thick]
%      (YOX) edge[-|,shorten <= 3pt, shorten >= 3pt, bend left] node[] {} (SWI)
%      (SWI) edge[->,shorten <= 3pt, shorten >= 3pt, bend left] node[] {} (YOX)
%      (SWI) edge[->,shorten <= 3pt, shorten >= 3pt] node[] {} (HCM)
%      (HCM) edge[->,shorten <= 3pt, shorten >= 3pt] node[] {} (NDD)
%      (NDD) edge[->,shorten <= 3pt, shorten >= 3pt, bend left] node[] {} (SWI)
%      ;

%      \end{tikzpicture}
%    \end{center}
%  \caption{A hypothetical regulatory network for the data in Figure~\ref{fig:dataset} where $M_{\text{SWI4}}$ is multiplication.}\label{fig:regnetwork}
%\end{figure}

The regulatory network $\mathbf{RN}_W$ shown in Figure~\ref{fig:regnetwork} is perhaps the simplest representative of the family of wavepool models proposed by the Haase lab~\cite{orlando2008global} for the metabolic cycle in \textit{S. cerevisiae}.
Our goal is to identify if, for a particular identification of the nodes $\setof{1,2,3,4}$  with the genes $\setof{ \text{SWI4, HCM1, NDD1, YOX1}}$, a DSGRN model of this form is consistent with the time series data shown in Figure~\ref{fig:dataset}(left), and, if so, under what ranges of parameter values.

Applying DSGRN database code to $\mathbf{RN}_W$ produces a parameter graph $\mathsf{PG}$ with 1080 nodes.
As explained in Section~\ref{section:DSGRN}, the phase space of this network is $(0,\infty)^4$ and the parameter space is a subset of $(0,\infty)^{19}$. 
The nodes correspond to 1080 distinct regions of parameter space which in turn give rise to 1080 distinct classes of state transition graphs which may arise from the regulatory network of Figure~\ref{fig:regnetwork}.
For each node we may present the associated non-empty connected region of parameter space as the solution set of a system of polynomial inequalities. 
For each point $z$ in this set, the  $DSGRN(\mathbf{RN}_W,z)$ system of trajectories gives rise to the same associated search graph.

\subsubsection{Invalidating a model.}
As a simple test we begin by considering a model that can be ruled out based on known biological interactions.
Consider the regulatory network $\mathbf{RN}_W$ where $1\leftrightarrow \text{NDD1}$, $2\leftrightarrow \text{HCM1}$, $3\leftrightarrow \text{SWI4}$, and $4\leftrightarrow \text{YOX1}$.
Applying our pattern matching methodology to the search graphs which arise for each of the 1080 parameter nodes corresponding to this instantiation of the regulatory network $\mathbf{RN}_W$ and the pattern graph of Figure~\ref{fig:patterngraph} we obtain no matches.
This indicates that no matter how parameters are chosen, the dynamical model cannot give rise to a solution trajectory exhibiting a behavior qualitatively similar to the collected experimental data of Figure~\ref{fig:dataset}.
Accordingly, we reject the proposed regulatory network model. %of Figure~\ref{fig:regnetwork_invalidate}.

\subsubsection{Parameter learning.}
We now turn to an accepted version of the wavepool regulatory network model $\mathbf{RN}_W$ where $1\leftrightarrow \text{SWI4}$, $2\leftrightarrow \text{HCM1}$, $3\leftrightarrow \text{NDD1}$, and $4\leftrightarrow \text{YOX1}$.
For this network we expect to find matches (in fact,  failure to find any matches would probably suggest that the DSGRN model was inappropriate).

Applying the pattern matching methodology to the search graphs which arise for each of the 1080 parameter nodes corresponding to this revised instantiation of the regulatory network $\mathbf{RN}_W$ and the pattern graph of Figure~\ref{fig:patterngraph} results in matches for 22  parameter nodes.  
By Theorem~\ref{thm:NoFalseNegatives}, for any parameter $z$ belonging to any of the other 1058 parameter nodes, the  $DSGRN(\mathbf{RN}_W,z)$ system of trajectories does not contain any trajectory passing only through domains and walls which exhibits a sequence of extrema matching a plausible total order of the experimentally observed extrema in the data.
Hence, our analysis has dramatically reduced uncertainty about relationships between the underlying parameters.

Furthermore, we can explicitly describe the regions of parameter space that correspond to these 22 matching parameter nodes.
For example, for one such parameter node the associated parameter region in $(0,\infty)^{19}$ is given by the inequalities 
\begin{eqnarray*}
& 0 < l_1 l_2 < \gamma_1\theta_3 < u_1l_2 < \gamma_1\theta_5 < l_1u_2 < u_1u_2 \\
& 0 < l_3 < \gamma_2 \theta_4 < u_3 \\
& 0 < l_4 < \gamma_3 \theta_1 < u_4 \\
& 0 < l_5 < \gamma_4 \theta_2 < u_5
\end{eqnarray*}
where 
\begin{eqnarray*}
l_1 := l((NDD1,SWI4)) & u_1 := u((NDD1,SWI4)) & \theta_1 := \Theta((NDD1,SWI4)) \\
l_2 := l((YOX1,SWI4)) & u_2 := u((YOX1,SWI4)) & \theta_2 := \Theta((YOX1,SWI4)) \\
l_3 := l((SWI4,HCM1)) & u_3 := u((SWI4,HCM1)) & \theta_3 := \Theta((SWI4,HCM1)) \\
l_4 := l((HCM1,NDD1)) & u_4 := u((HCM1,NDD1)) & \theta_4 := \Theta((HCM1,NDD1)) \\
l_5 := l((SWI4,YOX1)) & u_5 := u((SWI4,YOX1)) & \theta_5 := \Theta((SWI4,YOX1)) \\
\gamma_1 := \gamma(SWI4) & \gamma_2 := \gamma(HCM1) \\
\gamma_3 := \gamma(NDD1) & \gamma_4 := \gamma(YOX1).
\end{eqnarray*}
A complete listing of such regions is available in supplementary material \cite{SupplementalCodeRepository}.

\begin{figure}[h!]
    \begin{center}
	\includegraphics[width=1.75in]{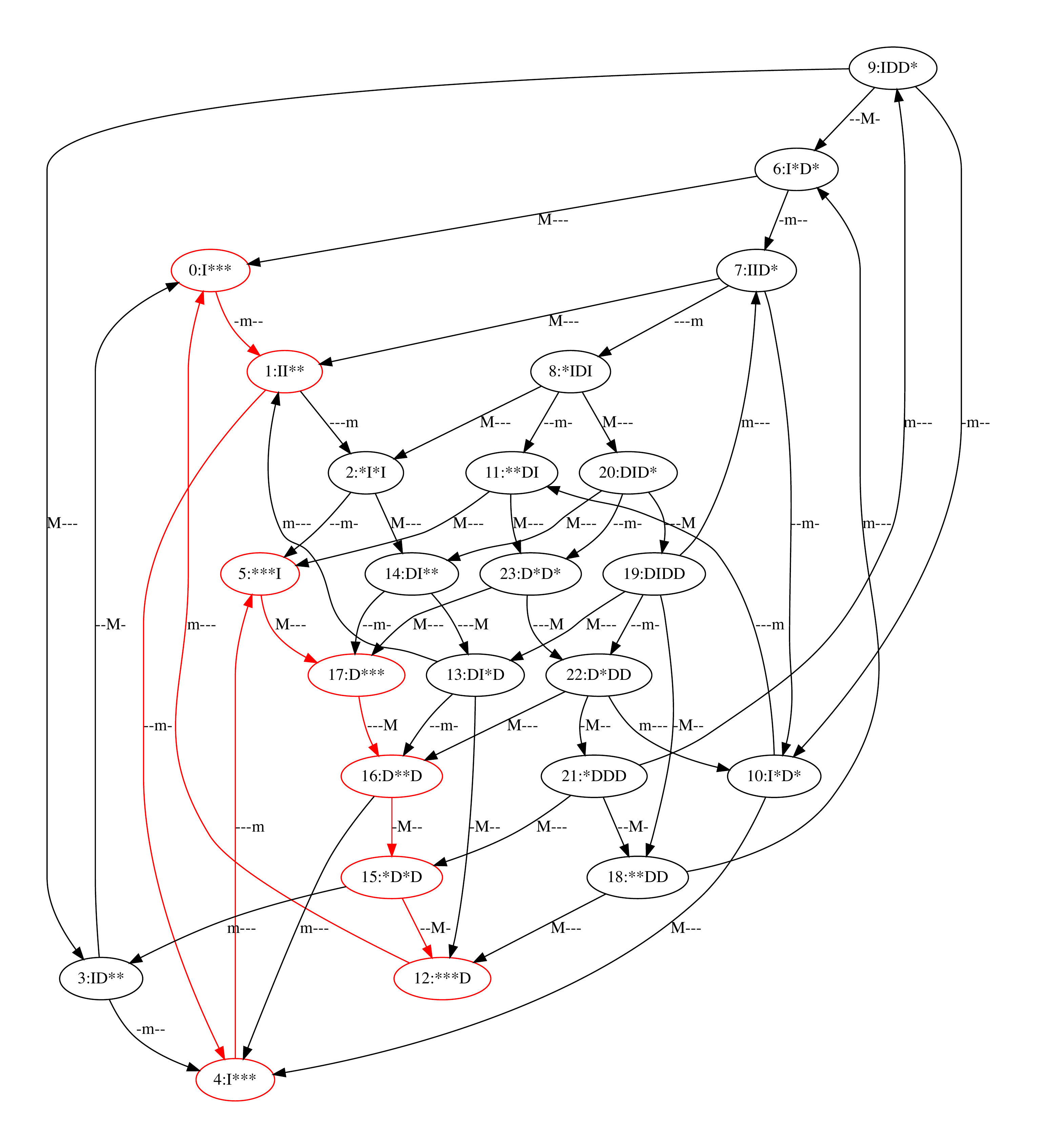}
	\includegraphics[width=1.75in]{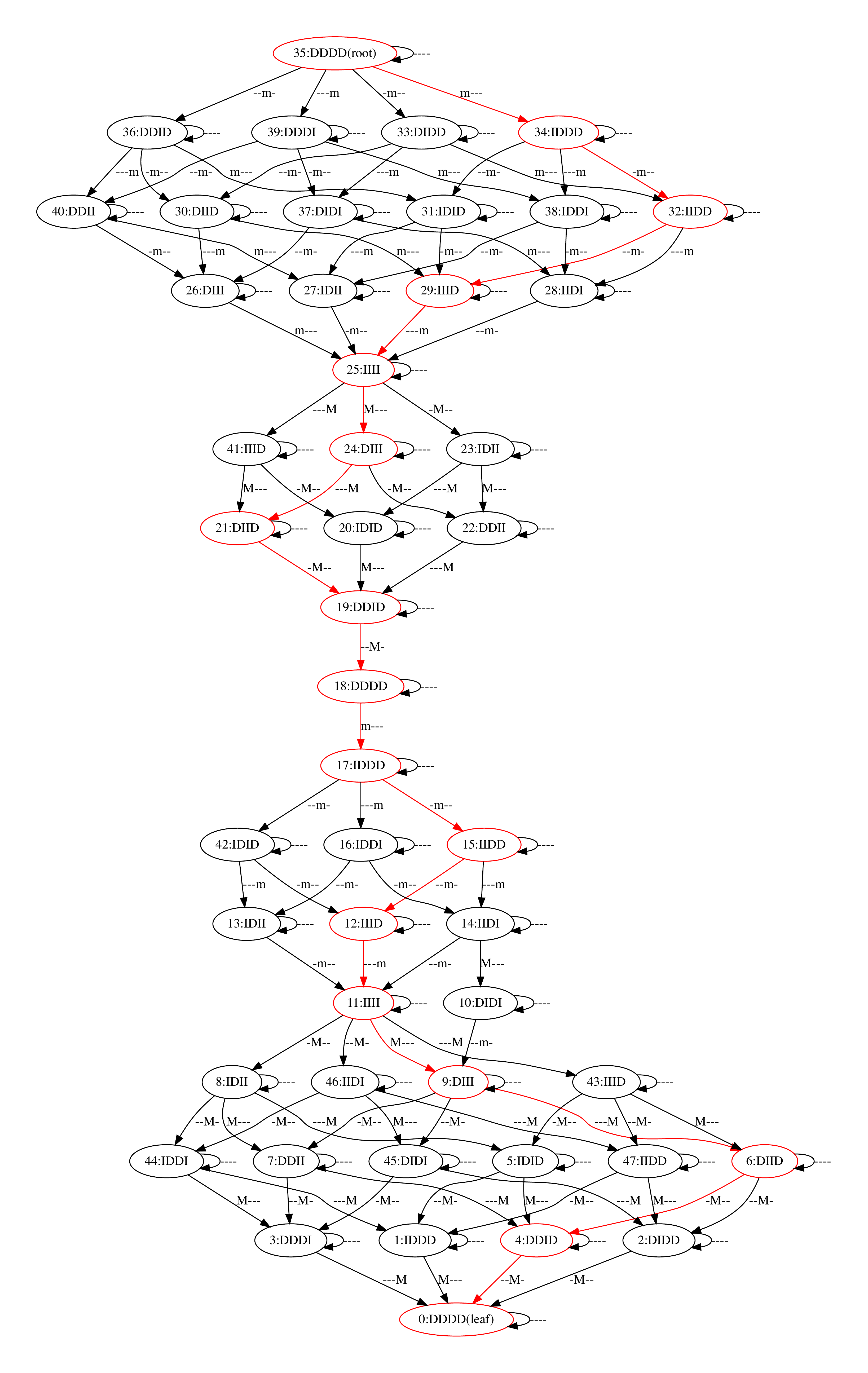}
	\includegraphics[width=1.75in]{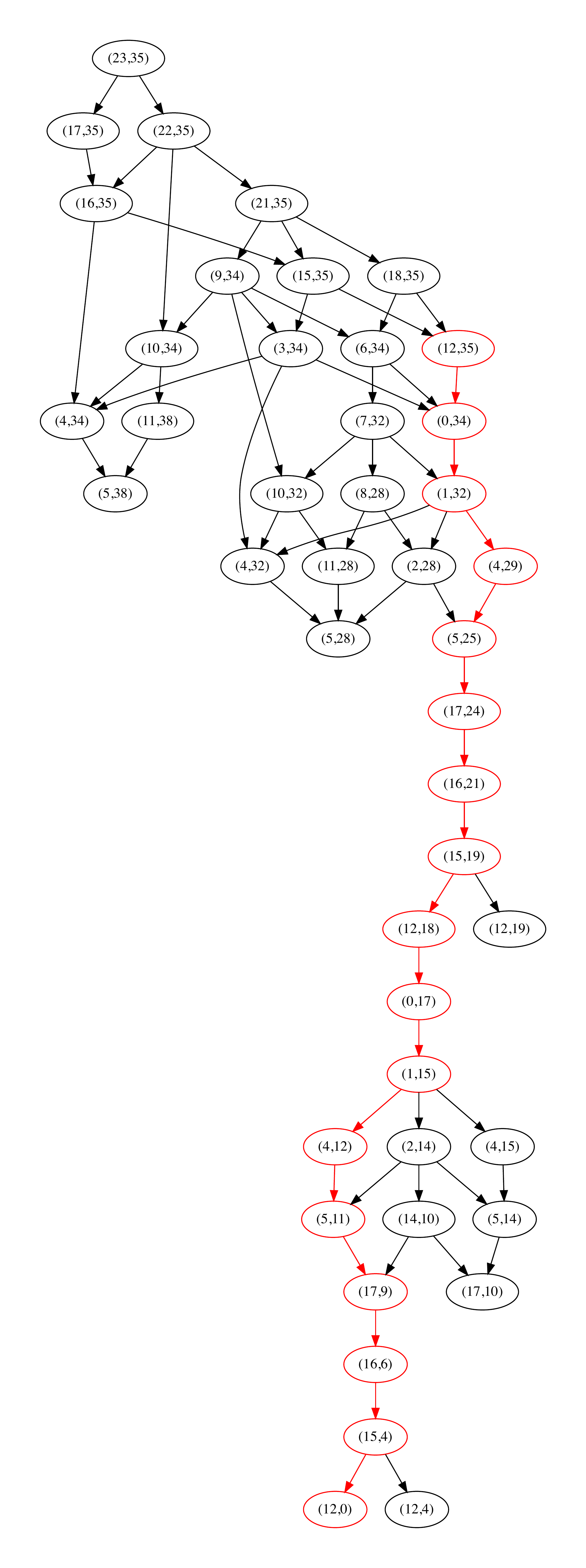}
	\caption{Left and Middle: Matching paths in search graph and pattern graph by Algorithm~\ref{alg:alignment}. Right: The corresponding path the algorithm found in the alignment graph.}
	\label{fig:matching_path}
	\end{center}
\end{figure}

A  pair of matching paths between the pattern graph and the search graph corresponding to this parameter region is shown in Figure~\ref{fig:matching_path}.

\section{Concluding Remarks}
\label{section:conclusion}
%!TEX root = ../PatternMatchPaper.tex
We presented a general method capable of rejecting models that cannot match coarse data generated by an experimentally measured time series. 
Our assumptions are very general; we expect that the time series is subject to substantial experimental error and therefore we only assume partial knowledge of the order of extrema of the components of the time series. This information is encoded in a poset of extrema, which we represent as a labeled directed acyclic graph called a {\em pattern graph}. 

Coming from the modeling side, we start with a concept of {\em system of trajectories}. Such a system can  be produced  by decomposition of the phase space into disjoint domains in which all trajectories are monotone, and on each boundary between domains, at most one component can attain an extremum. 
Existence of such a decomposition  allows extraction of the extremal behavior and its encoding into a {\em search graph}. 
On this level of generality we show that the problem of matching labeled paths between pattern graph and search graph can be solved in polynomial time. 

We discuss the applicability of our approach in two directions. First, we provide an example of a class of models which can be used to construct search graphs. 
Second, we apply our method to expression time series data from cell cycle in yeast. We show how our method can be used to learn  parameter regimes consistent with 
the experimental measurement by rejecting parameter regimes where the dynamics does not align with the data.

In order to ensure our results may be reproduced we adhere to the following recipe: (1) we release our code under an open-source license, (2) we host our code on a publicly available site using version-control (i.e. history tracking), (3) we give the version numbers of the code used to produce the result, (4) we provide instructions for installing and running the code, and (5) we produce \emph{digital object identifiers} (DOIs) of the versioned code for use in bibliographical entries.

The computer codes used to reproduce the results in this paper are stored in two code repositories. The first repository is the DSGRN project \cite{DSGRNRepository}. This is an open-source project which, as of writing, is hosted on the code-sharing website GitHub at \url{https://github.com/shaunharker/DSGRN}. The version utilized for this paper is 1.0.0. The second repository is the supplemental for this paper \cite{SupplementalCodeRepository} and houses the code (which relies on DSGRN) which is used to reproduce the above results. This again is open-source and is hosted at \url{https://github.com/shaunharker/2017-DSGRN-ModelRejection}. The version utilized for this paper is 1.0.0. The DOIs for these can be found in the references.

\section*{Acknowledgements}
The work of S. H. and K. M. was partially supported by grants NSF-DMS-1125174, 1248071, 1521771 and a DARPA contract HR0011-16-2-0033. T. G. was partially supported by  NSF grants DMS-1226213, DMS-1361240, DARPA D12AP200025 and NIH R01 grant 1R01AG040020-01.
B.C. was partially supported by DARPA D12AP200025 and USDA 2015-51106-23970.

\bibliographystyle{siamplain}
\bibliography{bibliography}
\end{document}